\documentclass[journal,twocolumn,twoside,10pt]{IEEEtran}

 \usepackage{amsmath,amssymb,amsthm,amsfonts,amsopn,cite,mathrsfs,todonotes}
\usepackage[latin1]{inputenc}
\usepackage[T1]{fontenc}
\usepackage{ifpdf}
\usepackage[english]{babel}
\usepackage[ruled]{algorithm2e}
\usepackage{dsfont,url}
\usepackage{tikz}
\usepackage{multirow}
\usetikzlibrary{patterns}
\usetikzlibrary{calc}
\usetikzlibrary{arrows.meta}

\newcounter{tempEquationCounter}
\newcounter{thisEquationNumber}

\theoremstyle{plain}
\newtheorem{theorem}{Theorem}

\theoremstyle{definition}
\theoremstyle{remark}
\theoremstyle{remark}
\newtheorem*{rem*}{Remark}

\newcommand{\bd}{\mathbf}   
\renewcommand{\mod}{\mathop{\operatorname{mod}}}  

\newcommand{\LtR}{\mathbf L^2(\mathbb{R})}

\newcommand{\LiR}{\mathbf L^1(\mathbb{R})}

\newcommand{\RR}{\mathbb R}

\newcommand{\NN}{\mathbb N}
\newcommand{\CC}{\mathbb C}

\def\sgn{\mathop{\operatorname{sgn}}}

\DeclareFontFamily{U}{mathx}{\hyphenchar\font45}
\DeclareFontShape{U}{mathx}{m}{n}{
      <5> <6> <7> <8> <9> <10>
      <10.95> <12> <14.4> <17.28> <20.74> <24.88>
      mathx10
      }{}
\DeclareSymbolFont{mathx}{U}{mathx}{m}{n}
\DeclareFontSubstitution{U}{mathx}{m}{n}
\DeclareMathAccent{\widecheck}{0}{mathx}{"71}
\DeclareMathAccent{\wideparen}{0}{mathx}{"75}
\newcommand{\grad}{\nabla}

\newcommand{\Parx}{\frac{\partial}{\partial x}}
\newcommand{\Parxt}{\frac{\partial^2}{\partial x^2}}
\newcommand{\Parxi}{\frac{\partial}{\partial \xi}}

\newcommand{\Pary}{\frac{\partial}{\partial y}}
\newcommand{\Paryt}{\frac{\partial^2}{\partial y^2}}

\usepackage{pgfplots,tikz}
\usepackage{subcaption}

\makeatletter
\renewcommand\section{\@startsection{section}{1}{\z@}{2ex plus 1.5ex minus 0.8ex}{1ex plus 1ex minus 0.5ex}{\normalfont\normalsize\centering\scshape}} 
\renewcommand\subsection{\@startsection{subsection}{2}{\z@}{2ex plus 1.5ex minus 0.8ex}{1ex plus 1ex minus 0.5ex}{\normalfont\normalsize\itshape}}
\renewcommand\@IEEEauthorblockconfadjspace{0.2em} 
\renewcommand\@IEEEauthorblockNtopspace{0.ex} 
\renewcommand\@IEEEauthorblockAtopspace{1ex}

\makeatother

\definecolor{tblblue}{rgb}{0.93,0.93,1.0}
\definecolor{tblred}{rgb}{1,0.93,0.93}
\definecolor{darkblue}{rgb}{0,0,0.5} 
\definecolor{darkgreen}{RGB}{20,120,43} 
\definecolor{darkred}{rgb}{0.8,0,0} 
\definecolor{lightblue}{RGB}{101,124,191}
\definecolor{skyblue}{RGB}{135,206,235}
\definecolor{gold}{RGB}{204,168,66}
\definecolor{strongblue}{RGB}{60,146,228}
\definecolor{lightgray}{gray}{0.5}
\definecolor{verylightgray}{RGB}{101,124,191}
\definecolor{mistyrose}{RGB}{238,213,210}
\definecolor{firebrick3}{RGB}{205,38,38}

\def\mwlet{\psi}
\def\magmwlet{M_{\mwlet}^s}
\def\magmwlettil{\tilde{M}_{\mwlet}^s}
\def\mygamma{\gamma}
\def\regamma{\operatorname{Re}(\gamma)}
\def\imgamma{\operatorname{Im}(\gamma)}

\allowdisplaybreaks[4]

\begin{document}
%
\title{Characterization of Analytic Wavelet Transforms and a New Phaseless
Reconstruction Algorithm
}
\author{Nicki Holighaus,
				G\"unther Koliander, 
        Zden\v{e}k Pr\r{u}\v{s}a, and Luis Daniel Abreu
\thanks{N.\ Holighaus*, G.\ Koliander, Z.\ Pr\r{u}\v{s}a, and L.\ D.\ Abreu are with the
Acoustics Research Institute, Austrian Academy of Sciences,
Wohllebengasse 12--14, 1040 Vienna, Austria, email:
\texttt{nicki.holighaus@oeaw.ac.at} (corresponding address),
\texttt{guenther.koliander@oeaw.ac.at},
\texttt{zdenek.prusa@oeaw.ac.at},
\texttt{labreu@kfs.oeaw.ac.at}}%
\thanks{Accompanying web page (sound examples, Matlab code, color figures)
\texttt{http://ltfat.github.io/notes/053}}%
\thanks{This work was supported by the Austrian Science Fund
(FWF): Y~551-N13, I~3067-N30, and P~31225-N32 and the Vienna Science and Technology Fund (WWTF): MA16-053.}%
\thanks{\textcopyright 2019 IEEE.  Personal use of this material is permitted.  Permission from IEEE must be obtained for all other uses, in any current or future media, including reprinting/republishing this material for advertising or promotional purposes, creating new collective works, for resale or redistribution to servers or lists, or reuse of any copyrighted component of this work in other works.}%
}

\markboth{}%
{
Holighaus, Koliander, Pr\r{u}\v{s}a, and Abreu: 
}
%


\IEEEaftertitletext{\vspace{-1\baselineskip}}
\maketitle

\begin{abstract}
We obtain a characterization of all wavelets leading to analytic wavelet
transforms (WT). The
characterization is obtained as a by-product of the theoretical foundations
of a new method for wavelet phase reconstruction from magnitude-only
coefficients. The cornerstone of our analysis is an expression of the
partial derivatives of the continuous WT, which results in
phase-magnitude relationships similar to the short-time Fourier transform
(STFT) setting and valid for the generalized family of Cauchy wavelets. We
show that the existence of such relations is equivalent to analyticity of
the WT up to a multiplicative weight and a scaling of the
mother wavelet. The implementation of the new phaseless reconstruction
method is considered in detail and compared to previous methods. It is shown
that the proposed method provides significant performance gains and a great
flexibility regarding accuracy versus complexity. Additionally, we discuss
the relation between scalogram reassignment operators and the wavelet
transform phase gradient and present an observation on the phase around
zeros of the WT.
\end{abstract}

\begin{IEEEkeywords}
gradient theorem, numerical integration, phase reconstruction,
short-time Fourier transform, wavelet transform, phase derivative, Cauchy-Riemann equations.
\end{IEEEkeywords}
\IEEEpeerreviewmaketitle

  \section{Introduction}
  
Time-frequency and time-scale representations are fundamental tools in many
areas of signal analysis and signal processing, ranging from medical data~%
\cite{nowak1999wavelet,chaplot2006classification}, damage or fault detection
in materials~\cite{marec2008damage} and machines~\cite{cusidocusido2008fault}%
, to image~\cite{chang2000adaptive,antonini1992image} and audio processing~%
\cite{yilmaz2004blind,chu2009environmental,necciari2018audlet}. Such representations are usually complex-valued and admit a natural
decomposition into magnitude and phase components, with the 
phase containing crucial information about the analyzed signal. However,
phase information is often discarded in favor of the magnitude information,
from which, supposedly, the desired information is more readily obtained. On
the other hand, whenever synthesis from the representation coefficients is
desired, the phase is crucial for the quality of the synthesized signal.

In recent years, the problems of phase retrieval and phaseless
reconstruction for time-frequency and time-scale dictionaries have attracted
considerable attention, leading to theoretical results for the feasibility
of phase retrieval in different contexts~\cite%
{mallatwaldspurger2017,balan2006signal,bandeira2014saving,alaifari2016stable,Alaifari2,Alaifari3}
and various algorithms~\cite%
{waldspurger2017,fienup1982phase,gerchberg1972practical,candes2013phaselift,griffin1984signal,shechtman2014gespar} (see also the survey \cite{shechtman2015phase}). 
Many
of these algorithms, e.g., \cite{griffin1984signal}, attempt to construct,
iteratively or directly, an appropriate phase to match the magnitude-only
representation coefficients, before finally performing a regular synthesis
step.

In many applications, in particular in audio signal processing, phaseless reconstruction
in time-frequency and time-scale dictionaries as proposed in, e.g., \cite{griffin1984signal,ltfatnote021}, 
is successfully applied. These applications include, but are not restricted to audio and speech 
synthesis~\cite{moulines1990pitch,wang2017tacotron,marafioti2019adversarial}, source 
separation~\cite{virtanen2007monaural,bach2006learning,bruna2015source}, as well as pitch and time-scale 
modifications~\cite{moulines1995non,laroche1999improved,pruuvsa2017phase}. In  these scenarios, 
either reconstruction from a phaseless representation is desired or the given phase has been 
invalidated in the course of processing and must be replaced.

Despite a relevant research activity regarding the phase behavior of
time-frequency and time-scale representations, with particular incidence in
the short-time Fourier transform (STFT)~\cite{auchfl12,po79,balazs2016pole}%
, applications that consider and/or modify the phase are relatively scarce.
Prime examples of analysis tools that successfully use phase information
are the so-called reassignment and synchrosqueezing methods~\cite%
{aufl95,damoto99,moto03,auger2013time,daubechies2011synchrosqueezed} that
deform the signal representation using a vector field obtained from the
phase gradient of the representation. 

For the STFT with Gaussian generator, the notion that
phase and magnitude carry equally important information has been made
precise by Portnoff~\cite{po79} and later by Auger and Flandrin~\cite%
{auchfl12}. They show that in this specific case, the phase gradient is
completely characterized by the gradient of the (logarithmically scaled)
magnitude.


\emph{Contributions: }
Our two main contributions are  a characterization of all analytic
wavelet transforms (WT) 
and  a new method for wavelet phase reconstruction from
magnitude-only coefficients. Here, the notion of analytic WTs is used in the sense that 
the WT is an analytic function of the upper half-plane (up to a signal-independent factor).
It is common to use the same terminology for wavelet transforms generated by analytic wavelets~\cite{liol10}, i.e., wavelets whose Fourier 
transform vanishes at negative frequencies. We discuss the connection between these two notions in detail.

The result characterizing all analytic wavelet
transforms shows that the most general wavelet leading to an analytic
transform has a Fourier transform given by
\begin{equation}
\xi^{\frac{\alpha -1}{2}}e^{-2\pi \gamma \xi }e^{i\beta \log \xi }
\label{analytic_wav}
\end{equation}%
for positive frequencies and we assume that the Fourier transform vanishes for negative frequencies.
The appearance of the $\beta$-dependent hyperbolic chirp in \eqref{analytic_wav} may be surprising since, to our knowledge, only the case $\beta =0$ 
has been associated with analytic functions in the literature (see \cite{dapa88}). 
It is commonly referred to as \emph{Cauchy wavelet} and historically 
associated with affine coherent states. 
However, the more general wavelets given by \eqref{analytic_wav} have been shown to minimize  a time-scale counterpart of Heisenberg uncertainty and are also known as ``Klauder wavelets'' \cite{fl98}.
The problem of
characterizing all analytic WTs has been considered and
partially solved in \cite{asbr09}, where it was also shown that
the Gaussian is essentially the only window leading to an analytic STFT.

We first discuss several aspects of the WT phase, loosely
following the structure of \cite{auchfl12}. Specifically, we express the
phase and log-magnitude derivatives in terms of
pointwise quotients of $3$ different WTs. The mother wavelets
used in these transforms can be derived directly from the mother wavelet of
the original WT. These expressions can be used, e.g.,
to estimate the local group delay and instantaneous scale (or frequency).
Furthermore, 
they impose conditions on the wavelets leading to analytic wavelet
transforms, and this is then used to show that 
the class of wavelets satisfying this condition are the generalized Cauchy
wavelets \eqref{analytic_wav}. The corresponding Cauchy Riemann (CR)
equations provide a relation between the phase gradient and the
(log-)magnitude gradient. 
Following \cite{balazs2016pole}, we also discuss the singular behavior
of the phase close to zeros of the WT and the implications of
our results for wavelet reassignment and ridge analysis.
Finally, we discuss the relation between analytic wavelets, i.e., wavelets  that can be extended to analytic functions on the upper half plane, and wavelets resulting in an analytic WT.

In the second part of the contribution, we implement a method for 
reconstruction from magnitude-only wavelet coefficients. More specifically, we
use a discrete approximation of the derived phase-magnitude relations to
modify the phase gradient heap integration algorithm~\cite%
{ltfatnote040,ltfatnote043,ltfatnote051}. The method is evaluated using
Cauchy wavelets of multiple orders, 
obtaining favorable results. In the evaluation, we further examine (uniform)
decimation and the number of scales at which the WT is
sampled, controlling the transform redundancy.

\emph{Structure of the Paper:} 
In Section~\ref{sec:PMCWT}, we present our expressions for the derivatives
of the the log-magnitude and the phase of the WT. These
expressions are used in Section~\ref{sec:pmrel} to obtain our result
characterizing the wavelets that lead to analytic WTs.
Furthermore, we prove direct relations between the log-magnitude and phase
derivatives. 
In Section~\ref{sec:CWasTF}, we provide a time-frequency interpretation of
the corresponding WT. Further applications to pole behavior
and scalogram reassignment of the expression provided in Section~\ref{sec:PMCWT} 
are presented in Sections~\ref{sec:polebeh} and \ref{sec:RArels}, 
respectively. 
In Section~\ref{sec:analyticvsAWT}, we discuss the relation between analytic wavelets and our analytic WT.
Finally, we apply the phase-magnitude relations to the
problem of phaseless reconstruction.
 In Section~\ref{sec:WPGHI}, we formally
introduce the discretization of our results and describe the phaseless
reconstruction algorithm; in Section~\ref{sec:experiments}, we conduct
experiments on real data and compare our method to previous approaches
toward phaseless reconstruction.

  \section{
	Log-Magnitude and Phase Derivatives of the Wavelet Transform
	}\label{sec:PMCWT}
  
  Fix a  function $\mwlet\in\LtR$ such that its Fourier transform $\widehat{\mwlet}$ vanishes almost everywhere on 
  $\RR^-$.  
  	The continuous WT (CWT) of a function (or signal) $s\in\LtR$ with respect to the \emph{mother wavelet} 
  $\mwlet$ is defined as
  \begin{equation}
    W_{\mwlet} s(x,y) = \langle s,\bd T_x \bd D_y \mwlet \rangle = \frac{1}{\sqrt{y}} \int_\RR s(t)\overline{\mwlet\left(\frac{t-x}{y}\right)}\, dt, 
		\label{eq:defwltransform}
  \end{equation}
  for all%
  \footnote{Although we restrict here to positive scales, there is no technical obstruction to allowing $y\in\RR\setminus\{0\}$.} 
  $x\in\RR$, $y\in\RR^+$. 
  Here, $\bd T_x$ and $\bd D_y$ denote the translation and dilation operators, respectively,  given by 
   $(\bd T_x s)(t) = s(t-x)$, and $(\bd D_y s)(t) = y^{-1/2} s(t/ y)$ for all $t\in \RR$.
  
  The CWT can be represented in terms of its magnitude $\magmwlet  := |W_\mwlet s| \geq 0$ and phase 
  $\phi^s_\mwlet := \arg (W_\mwlet s)\in\RR$ as usual. With this convention, $\log(W_\mwlet s) = \log(\magmwlet ) + i \phi^s_\mwlet$. This straightforward 
  relation is the basis of the following expressions for the partial derivatives of the log-magnitude and phase components, derived in  Appendix \ref{app:derivs}.
\begin{theorem}  
\label{thm:derexpr}
  Let  $\mwlet\in\LtR$ with $\widehat{\mwlet}(\xi)=0$ for $\xi< 0$ and assume that 
	$\mwlet$ is continuously differentiable with 
	$\mwlet', \bd T\mwlet'\in\LtR$,
	where $\bd T$, without subscript, denotes the time-weighting operator $(\bd T s)(t) = t s(t)$. 
  Then, for all $x\in \RR$ and $y\in \RR^+$ satisfying $\magmwlet (x,y)>0$, 
  \begin{equation}
    \nabla \log(\magmwlet )(x,y)
    = \begin{pmatrix}
    0 \\
    \frac{1}{2y}
    \end{pmatrix}
    - \operatorname{Re}   
    \begin{pmatrix}
    \frac{W_{\mwlet'} s(x,y) }{y W_\mwlet s(x,y)} \\[2mm] 
    \frac{W_{(\bd T \mwlet)'} s(x,y)}{y W_\mwlet s(x,y)} 
    \end{pmatrix}
	\label{eq:derlogmag}
  \end{equation}
	and
  \begin{equation}
    \nabla \phi^s_\mwlet(x,y)  = 
    - \operatorname{Im}
    \begin{pmatrix}
    \frac{W_{\mwlet'} s(x,y)}{y W_\mwlet s(x,y)} \\[2mm] 
    \frac{W_{(\bd T \mwlet)'} s(x,y)}{y W_\mwlet s(x,y)} 
    \end{pmatrix}.
	\label{eq:derphase}
  \end{equation}
	\end{theorem}

  The partial phase derivatives of the WT are often related to the local  instantaneous scale, as well as the 
  local group delay, of the analyzed signal, although at least the latter notion is not entirely clear 
  for WTs, see   Section~\ref{sec:RArels}. 
  The second order derivatives, which can be obtained similar to the first order derivatives in Appendix \ref{app:derivs},
  describe the variation of these quantities across phase space and thus proved useful as well, see~\cite{nelson2001cross,nelson2002instantaneous,auger2012making}. 
  %
  
	Formulas \eqref{eq:derlogmag}--\eqref{eq:derphase} can be used for the computation of the 
  partial derivatives by using efficient implementations of the WT.
	Furthermore, the possible accuracy of direct numerical differentiation of the 
  wavelet transform is limited if the WT can be computed only at
  certain positions 
  and not everywhere, e.g., in the presence of decimation in either coordinate.
  In contrast, even if a closed form expression for the derivative of the (time-weighted) wavelet is not available, numerical differentiation of the wavelet is not limited by this constraint.

\section{The Phase-Magnitude Relationship}\label{sec:pmrel}
  In general, the observations in Section~\ref{sec:PMCWT} do not yield a direct connection between the partial derivatives 
  of the log-magnitude and phase components. 
  However, if the WT is analytic, we can characterize the phase gradient by the log-magnitude gradient.
  Based on the Cauchy-Riemann (CR) equations,
  we can construct conditions on the mother wavelet $\mwlet$ such that the WT of any signal is an analytic function. 
  Similar to the analysis of analytic STFTs studied in \cite{asbr09} and the partial study of analytic wavelet transforms in the same contribution, we allow for an $(x,y)$-dependent factor $f(x,y)$ that is independent of the signal $s$.
  This factor can easily be accounted for when applying results from complex analysis to the transformed signal and leads to a significantly less restrictive class of  analyticity-inducing wavelets. 
  Furthermore, we want the class of analyticity-inducing wavelets to be invariant under the natural transforms associated with the WT, namely dilation and translation.
  Thus, we also allow for a constant dilation by $b \in \RR^+$ and a translation specified  by $a \in \RR$ in our analysis.
  \begin{theorem}\label{thm:analyticwlt}
    Let  $\mwlet\in\LtR$ with $\widehat{\mwlet}(\xi)=0$ for $\xi< 0$.
    There exist constants $a \in \RR$, $b \in \RR^+$ and a $\mathcal{C}^{\infty}$ function $f\colon \RR\times \RR^+\to \CC$ with $f(x,y)\neq 0$ such that 
    \begin{equation}
    \begin{split}
      h\colon \{z\in \CC: \operatorname{Im} (z) >0 \} &\to \CC  \\
       x+iy & \mapsto f(x,y)  W_\mwlet s(x-aby,by)
    \end{split}
    \label{eq:modwlettrans}
    \end{equation}
    is analytic for all $s\in \LtR$, 
    if and only if
    \begin{equation}
      \widehat{\mwlet}(\xi) 
      = c \xi^{\frac{\alpha-1}{2}} e^{- 2\pi\mygamma \xi} e^{i \beta \log  \xi}
      \label{eq:solutiondiffeq}
		\end{equation}
		for all $\xi >0$ and some constants $c\in \CC$, $\alpha > - 1$, $\beta \in \RR$, and  $\mygamma \in \CC$ with $\regamma >0$. 
  \end{theorem}
  A proof of the theorem is provided in Appendix \ref{app:analyticwlt}.
  The wavelets $\mwlet$ specified by \eqref{eq:solutiondiffeq}
  are known as Klauder wavelets \cite{fl98} and
  are a minor generalization of 
	Cauchy wavelets $\mwlet^{(\alpha)}$ \cite{dapa88}, which are recovered for the choice $\beta=0$ and $\mygamma=1$.
	Among other effects, modifying $\beta$ results in a proportional shift of the temporal concentration of the mother wavelet away from time zero.
  Furthermore, a change in $\mygamma$ results only in a scale change, dependent on $\regamma $, and a time shift, dependent on $\imgamma $. 
  Disregarding the constant factor $c\in\CC$ in \eqref{eq:solutiondiffeq}, we denote the generalized Cauchy wavelets in Theorem \ref{thm:analyticwlt} by $\mwlet^{(\alpha,\beta,\mygamma)}$
  or $\mwlet^{(\alpha,\beta)}$, if $\mygamma = 1$.

  If $\mwlet$ is given by \eqref{eq:solutiondiffeq}, a corresponding choice for $h$ being analytic is
  $f(x,y) = y^{-\frac{\alpha}{2}} e^{i \beta \log  y}$, $a=\imgamma $, and $b = 1/\regamma $, i.e.,
    \begin{equation}
       x+iy  \mapsto y^{-\frac{\alpha}{2}} e^{i \beta \log  y} W_\mwlet s\bigg(x- \frac{\imgamma }{\regamma } y,\frac{y}{\regamma } \bigg).
    \label{eq:genanalytic}
    \end{equation}
  We note that the wavelets are admissible only for $\alpha>1$.

  Theorem~\ref{thm:analyticwlt} can easily be modified to allow for wavelets $\mwlet$ where $\widehat{\mwlet}$ does not vanish for negative frequencies.
  In this case, $\widehat{\mwlet}$ must satisfy \eqref{eq:solutiondiffeq} for $\xi > 0$ and 
    \begin{equation}
      \widehat{\mwlet}(\xi) 
      = c_{\mathrm{n}} (-\xi)^{\frac{\alpha_{\mathrm{n}}-1}{2}} e^{ 2\pi\mygamma_{\mathrm{n}} \xi} e^{i \beta_{\mathrm{n}} \log (- \xi)}
    \end{equation}
  for $\xi<0$ and some constants $c_{\mathrm{n}}$, $\alpha_{\mathrm{n}}$, $\mygamma_{\mathrm{n}}$, and $\beta_{\mathrm{n}}$ satisfying the same constraints as $c$, $\alpha$, $\mygamma$, and $\beta$, respectively.

  Based on the CR equations, we obtain a phase-magnitude relation  for the WTs using $\mwlet^{(\alpha,\beta,\mygamma)}$. 
	\begin{theorem} 
	\label{thm:pmrel}
	Let $\mwlet$ be given by \eqref{eq:solutiondiffeq} with $c=1$.
	Then 
  \begin{equation}
  \label{eq:xDirPhaseDerivo}
	   \Parx  \phi_{\mwlet}^s
    =   \frac{\alpha}{2  y\regamma }
    - \frac{ \Pary  \log \big(\magmwlet\big) }{\regamma }
    + \frac{\imgamma\Parx  \log \big(\magmwlet\big)}{\regamma } 
	\end{equation}
	and
  \begin{align}
     \Pary  \phi_{\mwlet}^s 
    & = 
    \frac{\alpha\imgamma - 2 \beta}{2  y \regamma }   
    + \frac{ \lvert \mygamma\rvert^2  \Parx  \log  \big(\magmwlet \big)}{\regamma }
    \notag \\*
    & \quad 
    - \frac{\imgamma     \Pary  \log \big(\magmwlet\big)}{\regamma } .
  \label{eq:yDirPhaseDerivo}
  \end{align}
  For $\mygamma =1$, these relations simplify to
  \begin{equation}\label{eq:xDirPhaseDeriv}
     \Parx \phi^s_{\mwlet}(x,y) 
		= -  \Pary \log(\magmwlet )(x,y) + \frac{\alpha}{2y}
	\end{equation}
  and 
	\begin{equation}\label{eq:yDirPhaseDeriv}
    \Pary \phi^s_{\mwlet}(x,y) =  \Parx \log(\magmwlet )(x,y) - \frac{\beta}{y}.
	\end{equation}
	\end{theorem}
	A proof of the theorem is provided in Appendix \ref{app:PMrels}.
	As a simple consequence of Theorem~\ref{thm:pmrel}, we obtain for the second order derivatives in the case  $\mygamma =1$ 
    \begin{equation}
      \Parxt \log(\magmwlet )(x,y) +  \Paryt \log(\magmwlet )(x,y) 
      = - \frac{\alpha}{2  y^2}
    \end{equation}
    and
		\begin{equation}
      \Parxt \phi^s_{\mwlet}(x,y) +  \Paryt \phi^s_{\mwlet}(x,y)
			= \frac{\beta}{ y^2}.
		\end{equation}

    For Cauchy wavelets, it is known that the magnitude uniquely determines the phase up to a constant phase factor. Moreover, this statements even holds after  decimation in the scale component and in certain discretized settings \cite{mallatwaldspurger2017}.

	%

\section{The Cauchy Wavelet Transform as Time-Frequency Representation}\label{sec:CWasTF}

  If the mother wavelet $\mwlet$ is frequency-localized around frequency $\xi_{\mathrm{b}}$, then we can interpret $W_\mwlet s(x,y)$ 
  as a time-frequency measurement at frequency $\xi = \xi_{\mathrm{b}}/y$. In the case of the  wavelets $\mwlet^{(\alpha,\beta)}$,
  this leads to a particularly convenient form of the phase-magnitude relationship. 
  Here, we consider the unique peak of 
  $\lvert\widehat{\mwlet^{(\alpha,\beta)}}\rvert$ as center frequency, i.e., $\xi_{\mathrm{b}} = \frac{\alpha-1}{4\pi}$. Additionally, instead of the unitary
  dilation $\bd D_y$, we consider the dilation $\tilde{\bd D}_y s (t) = y^{-1} s(t/y)$ to define 
  \[
   \tilde{W}_\mwlet s(x,\xi) = \langle s,\bd T_x \tilde{\bd D}_{\xi_{\mathrm{b}}/\xi} \mwlet \rangle = \sqrt{\frac{\xi}{\xi_{\mathrm{b}}}}W_\mwlet s(x,\xi_{\mathrm{b}}/\xi).
  \]
  Using the relations in Section \ref{sec:pmrel},
	it is easy to derive
  \begin{align}\label{eq:xDirPhaseDerivAlt}
    \Parx \tilde{\phi}^s_{\mwlet}(x,\xi) 
	& = \frac{\xi^2}{\xi_{\mathrm{b}}} \Parxi \log(\magmwlettil)(x,\xi) + \frac{\alpha-1}{2\xi_{\mathrm{b}}}\xi\notag \\
	& = \frac{4\pi \xi^2}{\alpha-1}\Parxi \log(\magmwlettil)(x,\xi) + 2\pi\xi
  \end{align}
  and	
  \begin{align}\label{eq:xiDirPhaseDerivAlt}
     \Parxi \tilde{\phi}^s_{\mwlet}(x,\xi) & = -\frac{\xi_{\mathrm{b}}}{\xi^2} \left(\Parx \log(\magmwlettil)(x,\xi) - \frac{\beta\xi}{\xi_{\mathrm{b}}}\right)\notag \\
     & = -\frac{\alpha-1}{4\pi \xi^2} \Parx \log(\magmwlettil)(x,\xi)+\frac{\beta}{\xi},
  \end{align}
  where $\tilde{M}_\mwlet$ and $\tilde{\phi}_\mwlet$ denote the magnitude and phase of $\tilde{W}_\mwlet$, respectively.
  Interestingly, for standard Cauchy wavelets, i.e., $\beta=0$, these formulas coincide with the phase-magnitude relations for the STFT with a dilated Gaussian (see \cite[Section III]{ltfatnote040}) up to the simple
  change that the constant time-frequency ratio $\lambda$ is replaced by the
  frequency-dependent term $\frac{\alpha-1}{4\pi \xi^2}$. The above form 
  \eqref{eq:xDirPhaseDerivAlt}--\eqref{eq:xiDirPhaseDerivAlt} will enable us to adapt the 
  phase reconstruction method presented in \cite{ltfatnote051} to the WT more easily in Section~\ref{sec:WPGHI} below.
  
  The additive term $\beta/\xi$ in \eqref{eq:xiDirPhaseDerivAlt} compensates for the fact that $\mwlet^{(\alpha,\beta)}$ is time-localized around $\frac{2\beta}{\alpha-1}$. 
  In other words, the frequency (or scale) bands in $\tilde{W}_{\mwlet^{(\alpha,\beta)}}$ (or $W_{\mwlet^{(\alpha,\beta)}}$) are not time-aligned. Time-alignment can be restored by
  choosing $\imgamma = \frac{2\beta}{\alpha-1}$, thus removing the additive term at the cost of introducing directional derivatives of $\log(\magmwlettil)$ in the expression of the phase gradient:
  \begin{align}
    \Parx \tilde{\phi}^s_{\mwlet}(x,\xi) 
	& = \grad_{d_1(\xi)}\log(\magmwlettil)(x,\xi) + 2\pi\xi,	
  \end{align}
  with $d_1(\xi) = \left(\frac{2\beta}{\alpha-1},\frac{4\pi \xi^2}{\alpha-1}\right)$ and	
  \begin{align}
     \Parxi \tilde{\phi}^s_{\mwlet}(x,\xi) & = -\grad_{d_2(\xi)}\log(\magmwlettil)(x,\xi),
  \end{align}
  with $d_2(\xi) = \left(\frac{\alpha-1}{4\pi \xi^2}+ \frac{\beta^2}{(\alpha-1)\pi \xi^2},\frac{2\beta}{\alpha-1}\right)$.

\section{Derivatives and Analyticity of the Continuous Wavelet Transform---Further Observations}
\subsection{The Phase Around Zeros of the Wavelet Transform}\label{sec:polebeh}
  For the  STFT with  Gaussian window, 	
  it was remarked by Auger et al.~\cite{auchfl12}, that the phase has characteristic poles where the STFT is zero and that this fact can be derived from the analyticity of the Bargman transform. In 
  \cite{balazs2016pole}, the characteristic pole behavior was proven under weaker conditions, as long as the STFT is smooth enough, more specifically $\mathcal C^2$ or $\mathcal C^3$. 
  
  In fact, the techniques used therein apply to any complex-valued function of two real variables, 
  as long as its higher-order partial derivatives are continuous. In the case of the 
  WT, this can be ensured by selecting a sufficiently smooth and decaying 
  mother wavelet $\mwlet$. 
	In particular, if the $k$-th derivative of $\mwlet$ weighted by $t^l$ is square integrable, i.e., $\bd T^{l}\big(\frac{\partial^k}{\partial t^k}\mwlet \big) \in \bd L^2(\RR)$,
  for all $l,k\in\{0,\ldots,K\}$, then $W_\mwlet s\in \mathcal C^K(\RR\times\RR^+,\CC)$, for all
  $s\in\bd L^2(\RR)$, cf.\ Appendix~\ref{app:WTdiff}. 
	This implies the following result.
  
  \begin{theorem}\label{thm:wvltphase}
    Let $\mwlet\in\LtR$ and assume that $\bd T^{l}\big(\frac{\partial^k}{\partial t^k}\mwlet \big) \in \bd L^2(\RR)$,
    for all $l,k\in\{0,\ldots,2\}$. 
    If 
    $
     W_\mwlet s(x_0,y_0) = 0
    $
    and the Jacobian determinant $J_{x_0,y_0} :=det(D(W_\mwlet s))(x_0,y_0)$ of $W_\mwlet s$ at $(x_0,y_0)$ is nonzero, then for $\varepsilon > 0$ converging to $0$
    \begin{align}
      \sgn(J_{x_0,y_0}) \lim_{\varepsilon \to 0} \Parx \phi_\mwlet^s(x_0,y_0+ \varepsilon) & = \infty,
      \\
      \sgn(J_{x_0,y_0}) \lim_{\varepsilon \to 0} \Parx \phi_\mwlet^s(x_0,y_0 - \varepsilon) & = -\infty
    \end{align}
    and
    \begin{align}
      \sgn(J_{x_0,y_0}) \lim_{\varepsilon \to 0} \Pary \phi_\mwlet^s(x_0+ \varepsilon,y_0) & = \infty,
      \\
      \sgn(J_{x_0,y_0}) \lim_{\varepsilon \to 0} \Pary \phi_\mwlet^s(x_0 - \varepsilon,y_0) & = -\infty.
    \end{align}
    If even $\bd T^{l}\big(\frac{\partial^k}{\partial t^k}\mwlet \big) \in \bd L^2(\RR)$,
    for all $l,k \in \{0, \allowbreak \ldots, \allowbreak 3\}$, then the limits
    \begin{align}
    \lim_{\varepsilon \to 0} \Parx \phi_\mwlet^s(x_0+ \varepsilon,y_0)
    = \lim_{\varepsilon \to 0} \Parx \phi_\mwlet^s(x_0- \varepsilon,y_0)
    \end{align}
    and
    \begin{align}
    \lim_{\varepsilon \to 0} \Pary \phi_\mwlet^s(x_0,y_0+ \varepsilon)  
    = \lim_{\varepsilon \to 0} \Pary \phi_\mwlet^s(x_0,y_0- \varepsilon)
    \end{align}
    exist and are finite.
  \end{theorem}
  
  After noting that the assumptions imply that $W_\mwlet s\in\mathcal C^2(\RR\times\RR^+,\CC)$ (or $\mathcal C^3(\RR\times\RR^+,\CC)$), the proof of the above result is identical to the proofs of \cite[Theorem 4.7--4.9]{balazs2016pole}, which only rely on continuous differentiability locally.  
  
   \subsection{Scalogram Reassignment and Ridge Points}\label{sec:RArels}  
  
		Reassignment is a technique for sharpening time-frequency and time-scale representations \cite{aufl95}. 
		The reassignment map is a vector field that is used to deform the representation of choice and
    depends on the representation and the input signal.
		It has been derived using different 
    methods, e.g., as a constant phase deformation derived from group theoretical properties~\cite{damoto99,moto03} or by a center of
    gravity argument relying on the Wigner distribution~\cite{aufl95}. In the case of the spectrogram,
    both notions can be shown to lead to the same reassignment map.
		However, this is no longer true for the 
    scalogram, i.e., the squared modulus of the WT.
    
    The reassignment map given in~\cite{moto03} relies only on the phase gradient. It is defined as
    \begin{equation}\label{eq:ReassignPhase}
     (x,y) \mapsto \left(x+\frac{y^2\Pary \phi^s_{\mwlet}(x,y)}{\xi_{\mathrm{b}}}, \frac{\xi_{\mathrm{b}}}{\Parx \phi^s_{\mwlet}(x,y)}\right).
    \end{equation}
    Using \eqref{eq:derphase}, the map can be rewritten leading to 
    \begin{equation}
     (x,y) \mapsto \left(x-\frac{y\operatorname{Im}\left(\frac{W_{(\bd T \mwlet)'}s(x,y)}{W_{\mwlet}s(x,y)}\right)}{\xi_{\mathrm{b}}},-\frac{y\xi_{\mathrm{b}}}{\operatorname{Im}\left(\frac{W_{\mwlet'}s(x,y)}{W_\mwlet s(x,y)}\right)}\right).
    \label{eq:ReassignPhaseb}
    \end{equation}
    Similar to the expressions provided in \cite{aufl95}, \eqref{eq:ReassignPhaseb} can prove useful for the efficient calculation of reassigned scalograms.
    In particular, 
    a direct computation of the phase gradient from samples of the phase may be quite inaccurate.
		 If the phase-magnitude relations \eqref{eq:xDirPhaseDerivo} and \eqref{eq:yDirPhaseDerivo}
    are satisfied, then inserting these relations into the reassignment map \eqref{eq:ReassignPhase} 
    allows scalogram reassignment from the scalogram itself.

    In  \cite{aufl95}, a different reassignment map  is given by\footnote{Note that the equality $\mwlet = h(-\bullet)\exp(i\xi_{\mathrm{b}}(\bullet))$ converts 
    between the different WT conventions used.}
    \begin{equation}\label{eq:ReassignGrav}
     (x,y) \mapsto \left(x+y\operatorname{Re}\left(\frac{W_{\bd T \mwlet}s(x,y)}{W_{\mwlet}s(x,y)}\right),- \frac{y\xi_{\mathrm{b}}}{\operatorname{Im}\left(\frac{W_{\mwlet'}s(x,y)}{W_\mwlet s(x,y)}\right)}\right).
    \end{equation}
    The two reassignment maps coincide in the second coordinate, which can be considered a notion of local \emph{instantaneous scale}. 
		However, in the first
    coordinate, which is often considered an estimate of the local \emph{group delay}, they are quite different. 
		In particular, there seems to be no connection between the first coordinate of \eqref{eq:ReassignGrav} and partial derivatives of the WT $W_\mwlet s$ in general. 
    However, assuming $\bd T \mwlet = d_0\mwlet + d_1\mwlet'$ for some constants $d_0,d_1\in\CC$, enables an  expression of the first
    coordinate in \eqref{eq:ReassignGrav} as a linear combination of the partial derivatives of $\log (\magmwlet )$ and $\phi^s_{\mwlet}$. 
		A mother wavelet satisfying this differential equation 
    is, e.g.,  the \emph{Gabor wavelet} $\mwlet_{\mathrm{G}}(t) = e^{-t^2/2+i\xi_{\mathrm{b}} t}$, with $\bd T \mwlet_{\mathrm{G}} = i\xi_{\mathrm{b}} \mwlet_{\mathrm{G}} - (\mwlet_{\mathrm{G}})'$. 
		In this case, we can rewrite
    \begin{equation}
		\begin{split}
     \operatorname{Re}\left(\frac{W_{\bd T \mwlet_{\mathrm{G}}}s(x,y)}{W_{\mwlet_{\mathrm{G}}}s(x,y)}\right) 
		& = -\operatorname{Re}\left(\frac{W_{(\mwlet_{\mathrm{G}})'}s(x,y)}{W_{\mwlet_{\mathrm{G}}}s(x,y)}\right) \\
		& = y\Parx \log (M_{\mwlet_{\mathrm{G}}}^s)(x,y)
		\end{split}
    \end{equation}
    and therefore \eqref{eq:ReassignGrav} becomes
    \begin{equation}\label{eq:ReassignMorlet}
     (x,y) \mapsto \left(x+y^2\Parx \log (M_{\mwlet_{\mathrm{G}}}^s)(x,y), \frac{\xi_{\mathrm{b}}}{\Parx \phi^s_{\mwlet_{\mathrm{G}}}(x,y)}\right).
    \end{equation}   
    Note that the Gabor wavelet does not vanish at negative frequencies. However, the proof of Theorem \ref{thm:derexpr} does not rely on this property, such that the derivation of \eqref{eq:ReassignMorlet} from \eqref{eq:ReassignGrav} remains valid.
    
    The expression of the first coordinate in \eqref{eq:ReassignMorlet} is not too surprising, since $\bd D_y \mwlet_{\mathrm{G}}$ 
    is simply a dilated, modulated Gaussian, and thus this expression could also be obtained using the reassignment
    operators and phase-magnitude relationship for the Gaussian STFT, 
    as observed in~\cite{auchfl12}.
    
    Similar to reassignment, wavelet ridge analysis~\cite{liol10,delprat1992asymptotic} attempts to identify a \emph{skeleton} of essential time-scale positions in the wavelet transform. 
    The notion of magnitude ridge points is defined as the points $(x,y)$, such that  $\Pary \log (y^{-\frac{1}{2}}M_\mwlet^s)(x,y)=0$ 
    and $\Paryt  \log (y^{-\frac{1}{2}} M_\mwlet^s)(x,y) < 0$. 
    Similarly, phase ridge points are defined as the points $(x,y)$, such that  $\Parx \phi^s_{\mwlet}(x,y) - \frac{\xi_{\mathrm{b}}}{y}=0$ and $\Pary \left(\Parx \phi^s_{\mwlet}(x,y) - \frac{\xi_{\mathrm{b}}}{y}\right) > 0$. 
    It is straightforward to verify that Theorems \ref{thm:analyticwlt} and 
    \ref{thm:pmrel} imply that the phase and magnitude ridge points coincide if and only if $\mwlet = c\mwlet^{(\alpha,\beta,\gamma)}$, with $c\neq 0$, $\alpha>-1$, $\beta\in\RR$ 
    and $\gamma\in\RR^+$.   
    
\subsection{Analytic Wavelets and the Analytic Wavelet Transform}\label{sec:analyticvsAWT}

To prevent confusion related to other works on WTs, we want to point out the
connection between our analytic WT and the WT using analytic wavelets  \cite{liol10,liol12}.
  A wavelet is called analytic if it vanishes almost everywhere on 
  $\RR^-$, i.e., $\widehat{\mwlet}(\xi)=0$ for $\xi< 0$. 
  The reason for this terminology is that these wavelets can be extended to an analytic function on the upper half-plane by the Paley-Wiener theorem.
  Furthermore, the WT  using an analytic wavelet $\mwlet$ and at a fixed scale $y_0$ also has the property that it can be extended to an analytic function on the upper half-plane in a complex variable $w$, i.e., the function $x\mapsto W_{\mwlet} s(x,y_0)$ can be extended to an analytic function 
    \begin{equation}
      w\mapsto W^{(\textrm{a})}_{\mwlet} s(w,y_0).
      \label{eq:analyticwavelet}
    \end{equation}
  Now, for  $\mwlet=\mwlet^{(\alpha,\beta,\mygamma)}$, the function in \eqref{eq:genanalytic} is also analytic.
  The two functions in \eqref{eq:genanalytic} and \eqref{eq:analyticwavelet} coincide (up to a constant) for $y= \regamma y_0$ and $ x = w + \imgamma y_0$ for all $w\in \RR$.
  More specifically,
    \begin{align}
      & (\regamma y_0)^{-\frac{\alpha}{2}} e^{i \beta \log  (\regamma y_0)} W^{(\textrm{a})}_{\mwlet} s(w,y_0) 
      \notag \\ 
      & \qquad=
      h\big(w + \imgamma y_0+i(\regamma y_0)\big)
    \end{align}
  for all $w\in \RR$ and where both sides are analytic functions in $w$ on the upper half-plane.
  Thus, they have to coincide everywhere and we see that $h$ describes the analytic continuation $W^{(\textrm{a})}_{\mwlet} s(w,y_0)$ of an arbitrary scale up to some constants and shifts.
  In particular, for the case of Cauchy wavelets $\mwlet^{(\alpha)}$, we have
    \begin{equation}
       y_0^{-\frac{\alpha}{2}}  W^{(\textrm{a})}_{\mwlet} s(w,y_0) =
      h(w +i y_0)
    \end{equation}
  and the analytic continuation at a given scale only differs by a constant multiple from the analytic continuation at any other scale. 
  This \emph{equivalence} of all analytic continuations \eqref{eq:analyticwavelet} is  unique to the wavelets $\mwlet^{(\alpha,\beta,\mygamma)}$.
   
\section{Application---Phaseless Reconstruction for the Discrete Continuous Wavelet Transform}\label{sec:WPGHI}

In the following, we propose and evaluate a method for signal reconstruction from magnitude-only wavelet coefficients. 
More specifically, the proposed algorithm computes a phase estimate from the given magnitude-only coefficients. 
After combining the magnitude-only coefficients with the estimated phase, the wavelet transform must be inverted to obtain a time-domain signal. 
Here, any method that implements reconstruction from wavelet coefficients 
can be used.

Since
arbitrary dilations cannot be naturally transferred to the discrete domain, 
discrete implementations of the WT can be quite different from each other, 
see~\cite{mallat2008wavelet,rioul1992fast,unser1994fast,shensa1992discrete,ltfatnote038} and references therein. 
For illustrative purposes and to clarify notation, we shortly sketch an implementation of the discrete WT that follows~\cite{balazs2011theory,dogrhove13,schorkhuber2014matlab} closely, with some of the modifications introduced in \cite{necciari2018audlet}. In particular, we mimic the dilation operator by sampling the continuous frequency response 
of the mother wavelet $\mwlet\in\LtR\cap\LiR$ at the appropriate density. The following description reflects the implementation used in our experiments, see Section \ref{sec:experiments}. We use the terminology \emph{discrete continuous
WT (DCWT)} to distinguish this type of discrete WT from methods based on
a multiresolution analysis and wavelet bases~\cite{mallat2008wavelet,ltfatnote038}, that are commonly known as discrete wavelet transform (DWT).

\subsection{Discrete Continuous Wavelet Transform}
We will denote discretizations of continuous signals by brackets, e.g., the discretized signal $s_{\mathrm d}[l]\in \CC$ for $l\in \{1, \dots, L\}$ and some $L\in\NN$. 
In this discrete domain, the translation operator acts circularly, i.e., $s_{\mathrm d}[l-m]$ is 
interpreted as $s_{\mathrm d}[\mod(l-m,L)]$. Assuming the sampling rate $\xi_{\mathrm{s}}$, the wavelet at scale $y=\xi_{\mathrm{b}}/\xi$ is derived from the frequency response 
of the mother wavelet $\mwlet\in\LtR\cap\LiR$ as 
\[
 {\widehat{\mwlet_y}}[k] = {\widehat{\mwlet}}\left(\frac{y\xi_{\mathrm{s}} k}{L}\right) = {\widehat{\mwlet}}\left(\frac{\xi_{\mathrm{b}}\xi_{\mathrm{s}}}{L} \frac{k}{\xi}\right),
\]
for $k\in \{-\lfloor L/2 \rfloor, \dots, \lceil L/2\rceil -1\}$. 
Naturally, only a finite range 
of scales can be considered before the wavelet deteriorates either due to the 
sampling density being too coarse ($y$ large) or its bandwidth approaching $\xi_{\mathrm{s}}$ ($y$ small).
Hence, in order to cover the entire frequency range, we introduce an additional low-pass 
function in the style of~\cite[Section 3.1.2]{necciari2018audlet}. 

The entire wavelet system is characterized by the minimum scale 
$y_{\mathrm m}\in\RR^+$, 
the scale step\footnote{We choose the commonly used geometric spacing of center frequencies, but the proposed
phase reconstruction method remains valid for any center frequency spacing.} $2^{1/B}$, with $B\in\RR^+$, the number of scales $K\in\NN$, and
the decimation factor $a_{\mathrm d}\in\NN$, with $a_{\mathrm d}\vert L$.
The corresponding scaled and shifted wavelets are given as 
\begin{equation}
 \mwlet_{n,k} = \bd T_{na_{\mathrm d}} \mwlet_{2^{k/B} y_{\mathrm m}} 
\end{equation}
for $k \in  \{0, \dots, K-1\}$ and  $n \in  \{0, \dots, L/a_{\mathrm d}-1\}$.
A plateau function $P_{\mathrm{lp}}\in\CC^L$, centered at $0$, specifies the low-pass function as
\begin{equation}
  \widehat{\mwlet_{\mathrm{lp}}} =  a_{\mathrm d}^{-1} P_{\mathrm{lp}}\Psi_{\mathrm{lp}},
\end{equation}
where
\begin{equation}
 \Psi_{\mathrm{lp}} = \sqrt{\max(\Psi) - \Psi},\quad \Psi = \sum_{k=0}^{K-1} |\widehat{\mwlet_{0,k}}|^2.
\end{equation}

An analysis with the constructed system yields $ LK/a_{\mathrm d}$ com\-plex-valued coefficients for the wavelet scales and additional $L/a_{\mathrm d}$
real-valued coefficients for the low-pass function, for a total redundancy of $(2K+1)/a_{\mathrm d}$ when analyzing real-valued signals. With a slight abuse of 
terminology, we will from now on refer to the proportional quantity $K/a_{\mathrm d}$ as the redundancy. 

If $\mwlet$ is smooth and $a_{\mathrm d}, 1/B$ are small enough, then the results in~\cite{balazs2011theory} imply that the DCWT is invertible. 
Inversion can be achieved by interpreting the wavelet transform as a filter bank analysis and invoking the frame theory of uniform filter banks~\cite{bofehl98,cvetkovic1998oversampled,fickus2013finite,stroh1} to compute a 
dual filter bank synthesis. This can be done  either directly using dual filters $\widetilde{\mwlet}_{k}$, or iteratively~\cite{gr93,necciari2018audlet} using conjugate 
gradient iterations. The consideration of general uniform filter banks is necessary: It is not always possible to find a dual filter bank, which is required to achieve perfect reconstruction, with wavelet structure. Nonetheless, the dual filter bank shares the number of channels $K+1$ and the decimation factor $a$ of the wavelet analysis.  

  \subsection{Application to Phaseless Reconstruction}  
  For the phase-magnitude relations presented in Theorem~\ref{thm:pmrel} to hold, we have to assume a wavelet $\mwlet^{(\alpha,\beta,\mygamma)}$ as in 
  \eqref{eq:solutiondiffeq}. In particular, we will restrict to the case  $\mygamma =1$ for simplicity and drop the superscript $(\alpha,\beta)$ for notational convenience.
  For the WT and the phase-magnitude relations, we will use the convention introduced in Section~\ref{sec:CWasTF}. The generalization to the full class 
  of wavelets described by \eqref{eq:solutiondiffeq} is straightforward.
  
  Note that Theorem~\ref{thm:pmrel} provides only the phase-derivative and indeed reconstruction can at best be expected to be accurate up to a global phase factor.
  Furthermore, the reconstruction quality is expected to be worse for low-magnitude areas and thus the proposed algorithm only reconstructs the phase down to a certain magnitude-threshold. 
  Coefficients below that threshold are expected to have little effect on the synthesis and can thus be assigned a random phase. 
  As a consequence of the local, adaptive integration scheme, the reconstructed phase is in fact only expected to be consistent locally with changes by a constant phase factor between these local components. 
  On audio signals, such as the chosen corpus of test data, this change is not expected to have notable perceptual effects. 
  Nonetheless, it is visible in the phase difference (between original and reconstructed phase) in Fig.~\ref{fig:phasediffs}.
  At low redundancy, which is not well-suited for phase reconstruction in general, the phase distortion may become more severe (see the lower right corner in Fig.~\ref{fig:phasediffs}), sometimes leading to perceivable distortion.
  
  Assume that the continuous-time signal $s$ is approximately band- and time-limited on $[0,\xi_{\mathrm{s}})$ 
	and $[0,L/\xi_{\mathrm{s}})$, respectively. Then, with $s_{\mathrm d}[l] = s(l/\xi_{\mathrm{s}})$, for $l\in \{0, \dots, L-1\}$, 
  $a_{\mathrm d} = a\xi_{\mathrm{s}} \in\NN$, and $\xi_k = 2^{-k/B}\xi_{\mathrm{b}}/y_{\mathrm m}$, we obtain the approximation 
  \begin{equation}\label{eq:MagnitudeDisc}
    M_{s}[n,k] := |\langle s_{\mathrm d},\mwlet_{n,k} \rangle|  \approx 
    \xi_{\mathrm{s}}\magmwlettil (na,\xi_k).
  \end{equation} 
  Note that, after taking the logarithmic derivative of \eqref{eq:MagnitudeDisc},
  the normalization by $\xi_{\mathrm{s}}$ becomes irrelevant.
  
  Hence, we have 
  \begin{align}
    \Parx {\tilde{\phi}}^{s}_\mwlet(na,\xi_k)
    = &\ \frac{4\pi \xi_k^2}{\alpha-1}\Parxi \log(\magmwlettil )(na,\xi_k) + 2\pi \xi_k
    \notag \\*
    \approx \Delta^{\tilde{\phi},x,s}_\mwlet[n,k] 
    := &\ \frac{4\pi \xi_k^2}{\alpha-1}\Delta_k(\log(M_s))[n,k]+ 2\pi \xi_k,   
  \end{align}
  and
  \begin{align}
     \Parxi {\tilde{\phi}}^{s}_\mwlet(na,\xi_k) 
    = & -\frac{\alpha-1}{4\pi \xi_k^2}\Parx \log(\magmwlettil )(na,\xi_k)+\frac{\beta}{\xi_k}\notag \\*
    \approx \Delta^{\tilde{\phi},\xi,s}_\mwlet[n,k] 
    := & -\frac{\alpha-1}{4\pi \xi_k^2}\Delta_n(\log(M_s))[n,k]+\frac{\beta}{\xi_k}.   
  \end{align}
  Here, $\Delta_n$ and $\Delta_k$ are appropriate discrete differentiation schemes. 
  For $\Delta_n$, we can use centered differences, i.e.,
   \begin{equation}\label{eq:TDiffs}
     \Delta_n(M)[n,k] : = \frac{\xi_{\mathrm{s}}(M[n+1,k]-M[n-1,k])}{2a_{\mathrm d}}.
   \end{equation}
   The sampling step in the scale coordinate changes depends on $k$ and
   weighted centered differences can be used:
   \begin{align}\label{eq:SDiffs}
     \Delta_k(M)[n,k]
     & : = \frac{M[n,k+1]-M[n,k]}{2(\xi_{k+1}-\xi_k)}
     \notag \\
     & \quad 
     +\frac{M[n,k]-M[n,k-1]}{2(\xi_{k}-\xi_{k-1})}.
   \end{align}  
   
   Now, from $\Delta^{\tilde{\phi},x,s}_\mwlet$ and $\Delta^{\tilde{\phi},\xi,s}_\mwlet$, an estimate of the 
   phase of $\tilde{W}_\mwlet s$ at the sampling points $\{(na,\xi_k)\}_{n,k}$
   can be obtained using a quadrature rule considering the variable sampling intervals.
   That even simple $1$-dimensional trapezoidal quadrature provides satisfactory results  
   is illustrated by our experiments, see Section \ref{sec:experiments}. 
   
   The integration itself can be performed by a slightly 
   modified Phase Gradient
   Heap Integration (PGHI) algorithm~\cite{ltfatnote040,ltfatnote048,ltfatnote051},
   see Algorithm \ref{alg:wpghi}, using, e.g., the following integration rule on the set of neighbors $(n_{\mathrm{n}},k_{\mathrm{n}})\in \mathcal{N}_{n,k}:= \{(n\pm 1,k), (n,k\pm 1)\}$ of $(n,k)$
  \begin{align}is
    & (\tilde{\phi}^s_\psi)_{\mathrm{est}}[n_{\mathrm{n}},k_{\mathrm{n}}]
    \notag \\
    & = (\tilde{\phi}^s_\psi)_{\mathrm{est}}[n,k]
    +\frac{\xi_{k_{\mathrm{n}}}-\xi_k}{2}
    \left(\Delta^{\tilde{\phi},\xi,s}_\psi[n,k]+\Delta^{\tilde{\phi},\xi,s}_\psi[n_{\mathrm{n}},k_{\mathrm{n}}]\right)
    \notag \\
    & \quad 
    + \frac{a_{\mathrm d}(n_{\mathrm{n}} - n)}{2\xi_{\mathrm{s}}}
    \left(\Delta^{\tilde{\phi},x,s}_\psi[n,k]+\Delta^{\tilde{\phi},x,s}_\psi[n_{\mathrm{n}},k_{\mathrm{n}}]\right).
    \label{eq:twopointint}
  \end{align}
   When inserting \eqref{eq:TDiffs} and \eqref{eq:SDiffs} into \eqref{eq:twopointint}, the absolute scale of 
   the center frequencies $\xi_k$ and sampling rate $\xi_{\mathrm{s}}$ becomes unimportant and only their ratio
   enters the quadrature \eqref{eq:twopointint}. Hence, by considering relative frequencies $\xi_k/\xi_{\mathrm{s}}$,
   the algorithm is valid independent of the assumed sampling rate.

  \begin{algorithm}[htb]
    \LinesNumbered
    \caption{Wavelet Phase Gradient Heap Integration}
    \label{alg:wpghi}
    \KwIn{Magnitude $M_{s}$ of wavelet coefficients, estimates 
$\Delta^{\tilde{\phi},x,s}_\psi$
    and $\Delta^{\tilde{\phi},\xi,s}_\psi$ of the 
    partial phase derivatives, 
relative tolerance $\mathit{tol}$.}
    \KwOut{Phase estimate $(\tilde{\phi}^s_\psi)_{\text{est}}$.}    
    $\mathit{abstol} \leftarrow \mathit{tol}\cdot \max\left( M_{s}[n,k] \right)$\;
    Create set $\mathcal{I}=\left\{(n,k)  : 
M_{s}[n,k]>\mathit{abstol}
    \right\}$\;\label{alg:setI}
      Assign random values to $(\tilde{\phi}^s_\psi)_{\text{est}}(n,k)$ for $k\notin \mathcal{I}$\;
    Construct a self-sorting max \emph{heap} \cite{wi64} for $(n,k)$ 
pairs\;\label{alg:line}   
    \While{$\mathcal{I}$ \emph{is not} $\emptyset$  }{
        \If{heap \emph{is empty}}{
                Move 
$(n_m,k_m)=\smash{\mathop{\operatorname{arg~max}}\limits_{(n,k)\in\mathcal I}}
\left(M_{s}[n,k]\right)$ from $\mathcal I$ into the \emph{heap}\;
                $(\tilde{\phi}^s_\psi)_{\text{est}}(n_m,k_m) \leftarrow 0$\;               
 
        }
        \While{heap \emph{is not empty} }{
            $(n,k) \leftarrow$ remove the top of the
            \emph{heap}\;\label{alg:heapremove}
            
            \ForEach{$(n_{\mathrm{n}},k_{\mathrm{n}})$ \emph{in} $\mathcal{N}_{n,k} \cap \mathcal{I}$}{
                Compute $(\tilde{\phi}^s_\psi)_{\text{est}}(n_{\mathrm{n}},k_{\mathrm{n}})$ by
                means of \eqref{eq:twopointint}\;
                Move $(n_{\mathrm{n}},k_{\mathrm{n}})$ from  $\mathcal{I}$ into the \emph{heap}\;
            }
               
        }
    }  
    \bigskip
\end{algorithm}

Once the phase estimate $(\tilde{\phi}^s_\psi)_{\text{est}}$ has been computed, it is combined with the 
magnitude by $W_{s,\text{est}} := M_s e^{ i(\tilde{\phi}^s_\psi)_{\text{est}}}$. Subsequently, a time-domain 
signal can be synthesized as usual, e.g., using a dual filter bank. 

 \section{Experiments}\label{sec:experiments}

 To test and evaluate the proposed method, we performed two experiments, described and discussed
 below. Both experiments
 were run on the first $5$ seconds of all $70$ test signals from the Sound Quality Assessment Material recordings for subjective tests provided by the European Broadcasting Union (SQAM database)~\cite{sqam}.  
 For wavelet analysis and synthesis, we used the 
 filter bank methods in the open source Large Time-Frequency Analysis Toolbox (LTFAT~\cite{ltfatnote030}, \url{http://ltfat.github.io/}), where
 our implementation of \emph{Wavelet Phase Gradient
   Heap Integration (WPGHI)} is available by using the \texttt{'wavelet'} flag in 
 \texttt{filterbankconstphase}. A function to generate the wavelet filters and scripts for generating
 the individual experiments and figures are provided on the manuscript website 
 \url{http://ltfat.github.io/notes/053/}, where the resulting audio files for all experiment conditions can be 
 found as well.  
 Experimental conditions were restricted to classical Cauchy wavelets, i.e., $\beta = 0$ and $\mygamma =1$.
 
  Thus, the WT parameters used in the experiments are $(\alpha,a_{\mathrm d},K)$, where 
 $\alpha$ is the order of the Cauchy wavelet, $a_{\mathrm d}$ is the decimation step and $K$ is the number of frequency channels
 (or scales) used before adding the lowpass filter.
  As quantitative error measure, we employ (wavelet) spectral convergence~\cite{sturmel2011signal}, i.e., the relative mean squared error (in dB) between the wavelet coefficient magnitude of the target signal $s_{\mathrm t}$ and the proposed solution $s_{\mathrm p}$:
 \[
  SC(s_{\mathrm p},s_{\mathrm t}) = 20\log_{10} \frac{\|M_{s_{\mathrm p}}-M_{s_{\mathrm t}}\|}{\|M_{s_{\mathrm t}}\|}.
 \]
 It should be noted that the wavelet coefficient magnitude in the above formula was computed using the same parameter set $(\alpha,a_{\mathrm d},K)$ for which phaseless reconstruction was attempted.%
\footnote{Although spectral convergence is in some cases sensitive to parameter changes, preliminary tests showed that, for fixed $\alpha$, comparable spectral convergence is achieved with respect to representations with varying $a_{\mathrm d}$ and $K$. On the other hand, when $\alpha$ is changed as well, then the value of $SC(s_p,s_t)$ may change dramatically, such that comparing the results across different choices of $\alpha$ may be misleading.}
  
 \subsection{Experiment I---Comparison to Previous Methods}\label{ssec:exp1}
 
 To study the performance of the proposed algorithm in comparison with previous methods for phaseless recovery 
 from wavelet coefficients, we selected three settings of the WT parameters $(\alpha,a_{\mathrm d},K)$. 
For all settings, the channel center frequencies where geometrically spaced in $\frac{\xi_s}{20}\cdot [2^{-6},2^{3.3}]$. We considered the following tuples of parameters: $(30,5,100)$, $(300,12,240)$, and $(3000,20,400)$. Here, the ratio $K/a_{\mathrm d} = 20$ was fixed in all cases, but $a_{\mathrm d}$ and $K$ were adjusted to accommodate for bandwidth variations with changing $\alpha$. 
 
\begin{figure}[tbh]
\centering
\begin{subfigure}{0.47\textwidth}
\begin{tikzpicture}
\begin{axis}[
        legend columns=2,
        height=6.5cm,
        ymax=0, 
        xmin=0, xmax=71,
        tick label style={font=\small},
        width=\textwidth,
        grid=none,
        minor tick num=3,
        axis background/.style={fill=white},
        ylabel={Spectral Convergence},
        tick align=outside,
        legend entries={
				FBPGHI,
				WPGHI, 
        R-FGLIM,
        W-FGLIM
        },
		legend style={legend pos=north west, font=\footnotesize}
     ]   
     \addplot+[red!50,  mark=o, mark size=1pt, mark options={red!50}, line width=0.3pt] table[x expr=\thisrowno{0},y expr=(\thisrowno{2})] {compmethodsalpha30.dat};
     \addplot+[darkgreen!50,  mark=x, mark size=1.5pt, mark options={darkgreen!50}, line width=1pt] table[x expr=\thisrowno{0},y expr=\thisrowno{1}] {compmethodsalpha30.dat};
     \addplot+[darkblue, no markers, line width=0.3pt] table[x expr=\thisrowno{0},y expr=(\thisrowno{3})] {compmethodsalpha30.dat};
     \addplot+[black, mark=triangle*, mark size=1pt, mark options={black}, line width=1pt] table[x expr=\thisrowno{0},y expr=(\thisrowno{4})] {compmethodsalpha30.dat};
\end{axis}
\end{tikzpicture}
\end{subfigure}
\begin{subfigure}{0.47\textwidth}
\begin{tikzpicture}
\begin{axis}[
        legend columns=2,
        height=6.5cm,
        ymax=0, 
        xmin=0, xmax=71,
        tick label style={font=\small},
        grid=none,
        minor tick num=3,
        width=\textwidth,
        axis background/.style={fill=white},
        ylabel={Spectral Convergence},
        tick align=outside,
        legend entries={
				FBPGHI,
				WPGHI, 
        R-FGLIM,
        W-FGLIM
        },
		legend style={legend pos=north west, font=\footnotesize}
     ]   
     \addplot+[red!50,  mark=o, mark size=1pt, mark options={red!50}, line width=0.3pt] table[x expr=\thisrowno{0},y expr=(\thisrowno{2})] {compmethodsalpha300.dat};
     \addplot+[darkgreen!50,  mark=x, mark size=1.5pt, mark options={darkgreen!50}, line width=1pt] table[x expr=\thisrowno{0},y expr=\thisrowno{1}] {compmethodsalpha300.dat};
     \addplot+[darkblue, no markers, line width=0.3pt] table[x expr=\thisrowno{0},y expr=(\thisrowno{3})] {compmethodsalpha300.dat};
     \addplot+[black, mark=triangle*, mark size=1pt, mark options={black}, line width=1pt] table[x expr=\thisrowno{0},y expr=(\thisrowno{4})] {compmethodsalpha300.dat};
\end{axis}
\end{tikzpicture}
\end{subfigure}
\begin{subfigure}{0.47\textwidth}
\begin{tikzpicture}
\begin{axis}[
        legend columns=2,
        height=6.5cm,
        ymax=0, 
        xmin=0, xmax=71,
        tick label style={font=\small},
        width=\textwidth,
        grid=none,
        minor tick num=3,
        axis background/.style={fill=white},
        ylabel={Spectral Convergence},
        xlabel={Signal sorted by R-FGLIM performance},
        tick align=outside,
        legend entries={
				FBPGHI,
				WPGHI, 
        R-FGLIM,
        W-FGLIM
        },
		legend style={legend pos=north west, font=\footnotesize}
     ]   
     \addplot+[red!50,  mark=o, mark size=1pt, mark options={red!50}, line width=0.3pt] table[x expr=\thisrowno{0},y expr=(\thisrowno{2})] {compmethodsalpha3000.dat};
     \addplot+[darkgreen!50,  mark=x, mark size=1.5pt, mark options={darkgreen!50}, line width=1pt] table[x expr=\thisrowno{0},y expr=\thisrowno{1}] {compmethodsalpha3000.dat};
     \addplot+[darkblue, no markers, line width=0.3pt] table[x expr=\thisrowno{0},y expr=(\thisrowno{3})] {compmethodsalpha3000.dat};
     \addplot+[black, mark=triangle*, mark size=1pt, mark options={black}, line width=1pt] table[x expr=\thisrowno{0},y expr=(\thisrowno{4})] {compmethodsalpha3000.dat};
\end{axis}
\end{tikzpicture}
\end{subfigure}
\caption{70 signals for $\alpha = 30$ (top), $\alpha = 300$ (middle), and  $\alpha = 3000$ (bottom) sorted by R-FGLIM performance; the sort sequences are available on the manuscript webpage.\label{fig:perf30}}
\end{figure}
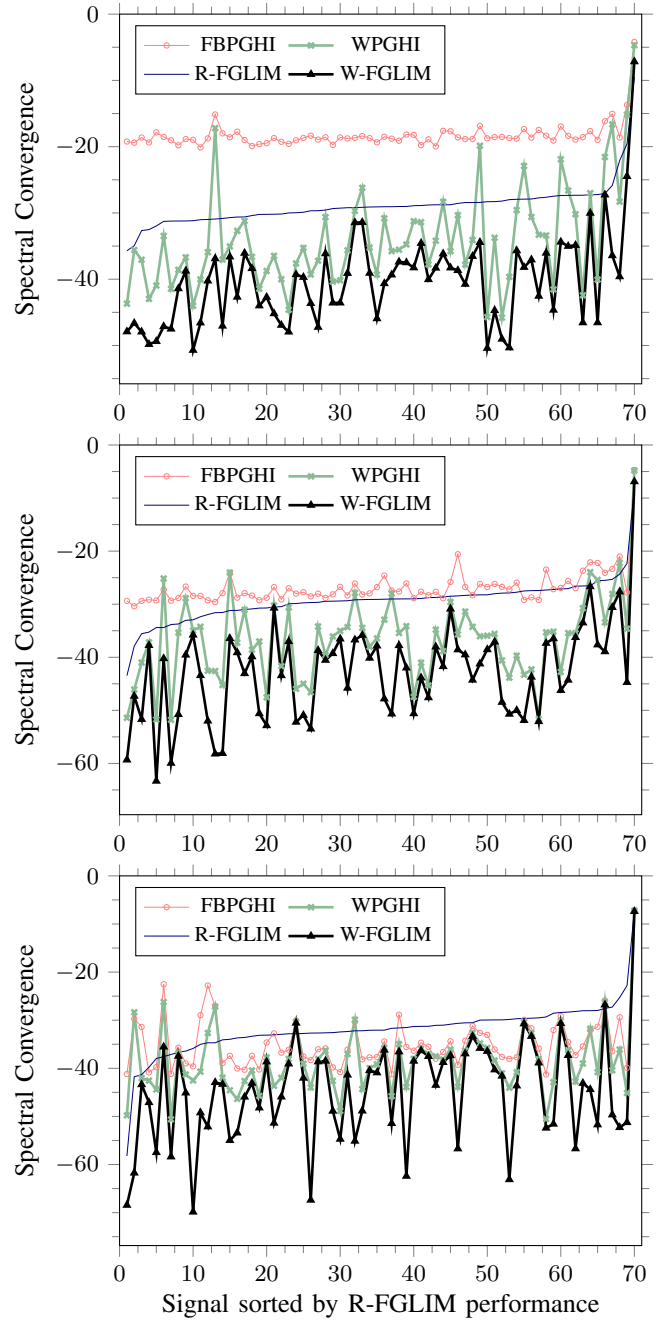

 The dimensionality of the considered audio data renders a systematic comparison to existing implementations of some established methods, e.g., \cite{candes2013phaselift,waldspurger2017} unfeasible, such that we
 resort to \emph{fast Griffin-Lim}~\cite{griffin1984signal,ltfatnote021} as baseline method.
 We compare four different methods: wavelet PGHI (WPGHI, \emph{proposed}), filter bank PGHI (FBPGHI, \cite{ltfatnote051}), fast Griffin-Lim with random initialization (R-FGLIM, \cite{ltfatnote021}) and fast Griffin-Lim initialized with the result of WPGHI (W-FGLIM, \emph{proposed}). 
Fast Griffin-Lim was restricted to at most $150$ iterations.%
\footnote{Although we are mainly interested in the reconstruction quality and not in computational performance, it is worth mentioning that the solutions of plain WPGHI and FBPGHI are computed in a small fraction of the time required for executing either R-FGLIM or W-FGLIM, even if, at the cost of reconstruction quality, the maximum number of iterations was significantly reduced.} 
  Spectral convergence of the four methods on all test signals is shown in Figure~\ref{fig:perf30} for the different parameter sets $(\alpha,a_{\mathrm d},K)$. 
  The means and standard deviation across all signals, for every method and parameters set are shown in Table \ref{tab:compMethods}. 
 
 It can be seen that on average, plain WPGHI (\emph{proposed}) outperforms both R-FGLIM and FBPGHI on all parameter sets, although FBPGHI approaches the other methods for larger values of $\alpha$.  
 Moreover, W-FGLIM (\emph{proposed}) shows significant improvements over either WPGHI or R-FGLIM. Looking at the individual signals more closely, we see in Figure~\ref{fig:perf30} that there are only very few cases in which R-FGLIM yields a better result than W-FGLIM.
For $\alpha = 3000$, all methods show comparable performance, with the exception of W-FGLIM, which still provides a clear advantage. 
 
 The figures and the computed standard deviations both suggest that methods that perform well on average are prone to larger performance fluctuation between individual signals. 
However, there is a small set of signals on which all methods perform badly, indicating that the fault is with the wavelet representation rather than the method applied for phaseless reconstruction. 
Notably, signal $65$ (corresponding to the rightmost signal in Figures~\ref{fig:perf30}--\ref{fig:redu2}) from the SQAM database, which, in the considered range, only contains a sustained, extremely low-pitched note is badly resolved by the employed WT and yields the worst spectral convergence of all signals, for any of the employed methods. 
 
 Informal listening showed that, for $\alpha = 3000$ and an untrained listener, the obtained reconstructions are mostly indistinguishable from the original signal.
 For $\alpha = 300$ and, more prominently, $\alpha = 30$, FBPGHI often introduces a characteristic pitch-shift, likely due to a wrongly estimated time-direction phase derivative, while R-FGLIM suffers from undesired frequency modulation artifacts; both types of distortion are most easily audible in simple signals, such as signal $1$ (sine wave) and $4$ (electronic gong) of the SQAM database and not present in the reconstructions provided by WPGHI. For $\alpha=30$ and some select cases, e.g., signals $16$ (clarinet) and $32$ (triangle),
audible distortions were present in the solutions by W-FGLIM, but not those by plain WPGHI, indicating that the observed improvement in terms of spectral convergence does not necessarily provide a perceptual improvement. To confirm our observations and to form their own opinion, the reader is invited to visit the manuscript webpage \url{http://ltfat.github.io/notes/053/}. 
 
 As a side note, the audio examples\footnote{Avalabile at \url{http://ltfat.github.io/notes/040/}} provided with~\cite{ltfatnote040} have clearly audible artifacts for signal $54$ (male German speech), for STFT-based PGHI and some competing algorithms. These artifacts are not present in any of the reconstructions we obtained using WPGHI, R-FGLIM, or W-FGLIM, for any considered parameter set, indicating that in some cases, usage of the WT may provide a genuine advantage over the STFT.
      
\begin{table}[thb] 
     \begin{center}\begin{tabular}{ |cc|cccc|} 
       \hline 
     $\alpha$ & method &  WPGHI &  FBPGHI & R-FGLIM & W-FGLIM \\ 
      \hline  
      \multirow{2}{*}{30} & mean & $ -33.9512 $ & $ -18.2916 $ & $ -28.9164 $ & $ -40.2932 $ \\ 
      & std & $   7.7632 $ & $   2.0358 $ & $   3.2641 $ & $   7.0526 $ \\ 
      \hline 
      \multirow{2}{*}{300} & mean & $ -36.6393 $ & $ -26.8642 $ & $ -29.3628 $ & $ -42.7550 $ \\  
      & std & $   7.9517 $ & $   3.4095 $ & $   4.0747 $ & $   9.0264 $ \\
      \hline 
      \multirow{2}{*}{3000} & mean & $ -38.5527 $ & $ -34.8240 $ & $ -32.0108 $ & $ -44.9260 $ \\ 
      & std & $   6.7961 $ & $   5.5137 $ & $   5.4821 $ & $  10.6349 $ \\ 
      \hline 
      \end{tabular}\end{center} 
      \caption{Means and standard deviation of spectral convergence for the considered methods and parameter sets.\vspace{-5mm}} 
      \label{tab:compMethods}
\end{table}      
 
  \subsection{Experiment II---Changing the Redundancy}
  
  In a second set of experiments, we investigate the influence of the redundancy $K/a_{\mathrm d}$ on the performance of the proposed 
  methods WPGHI and W-FGLIM. 
  For this purpose, we fixed an intermediate value for the order parameter, 
  setting $\alpha = 1000$ and consider redundancies $K/a_{\mathrm d} \in \{3,5,10,30\}$. 
  Here, however, the reconstruction quality does not only depend on the accuracy of WPGHI on the given magnitude coefficients, but also on the robustness of the synthesis by the dual system. 
  This robustness can be quantified by the so-called frame bound ratio of the respective wavelet system (for details  see \cite{ch16,da92,mallat2008wavelet}).
  
  In the ranges considered and for fixed $K/a_{\mathrm d}$, 
  the number $K$ of frequency channels has a larger influence on WPGHI performance than the decimation step $a_{\mathrm d}$, which was generally small. 
  On the other hand, the wavelet frame bound ratio deteriorates very quickly\footnote{Large frame bound ratios also decrease numerical stability, such that audio file generation from the obtained reconstructions is prone to clipping artifacts.} for too large decimation steps $a_{\mathrm d}$. Hence, the choice of wavelet parameters was a trade-off between the two factors with no clear optimal solution. After some preliminary testing, we fixed the following parameter sets $(\alpha,a_{\mathrm d},K)$: Low redundancy $(1000,30,90)$ (low), Medium redundancy $(1000,25,125)$ (medium), Medium high redundancy $(1000,18,180)$ (medhigh), High redundancy $(1000,10,300)$ (high).
  
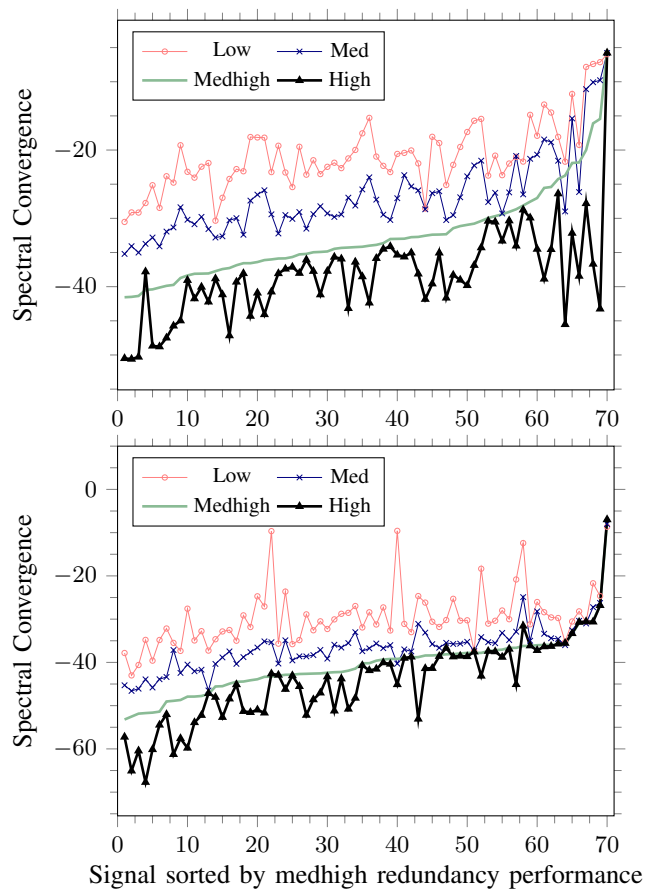
\begin{figure}[tb]
\centering
\begin{subfigure}{0.47\textwidth}
\begin{tikzpicture}
\begin{axis}[
        legend columns=2,
        height=6.5cm,
        xmin=0, xmax=71,
        tick label style={font=\small},
        width=0.96\textwidth,
        grid=none,
        minor tick num=3,
        axis background/.style={fill=white},
        ylabel={Spectral Convergence},
        tick align=outside,
        legend entries={
				Low,
        Med,
        Medhigh,
        High
        },
		legend style={legend pos=north west, font=\footnotesize}
     ]   
     \addplot+[red!50,  mark=o, mark size=1pt, mark options={red!50}, line width=0.3pt] table[x expr=\thisrowno{0},y expr=(\thisrowno{4})] {compredunwpghi.dat};
     \addplot+[darkblue,  mark=x, mark size=1.5pt, mark options={darkblue}, line width=0.3pt] table[x expr=\thisrowno{0},y expr=\thisrowno{3}] {compredunwpghi.dat};
     \addplot+[darkgreen!50, no markers, line width=1pt] table[x expr=\thisrowno{0},y expr=(\thisrowno{2})] {compredunwpghi.dat};
     \addplot+[black, mark=triangle*, mark size=1pt, mark options={black}, line width=1pt] table[x expr=\thisrowno{0},y expr=(\thisrowno{1})] {compredunwpghi.dat};
\end{axis}
\end{tikzpicture}
\end{subfigure}
%

\begin{subfigure}{0.47\textwidth}
\begin{tikzpicture}
\begin{axis}[
        ymax=10, 
        legend columns=2,
        height=6.5cm,
        xmin=0, xmax=71,
        tick label style={font=\small},
        width=0.96\textwidth,
        grid=none,
        minor tick num=3,
        axis background/.style={fill=white},
        ylabel={Spectral Convergence},
        xlabel={Signal sorted by medhigh redundancy performance},
        tick align=outside,
        legend entries={
				Low,
        Med,
        Medhigh,
        High
        },
		legend style={legend pos=north west, font=\footnotesize}
     ]   
     \addplot+[red!50,  mark=o, mark size=1pt, mark options={red!50}, line width=0.3pt] table[x expr=\thisrowno{0},y expr=(\thisrowno{4})] {compredunwpghifglim.dat};
     \addplot+[darkblue,  mark=x, mark size=1.5pt, mark options={darkblue}, line width=0.3pt] table[x expr=\thisrowno{0},y expr=\thisrowno{3}] {compredunwpghifglim.dat};
     \addplot+[darkgreen!50, no markers, line width=1pt] table[x expr=\thisrowno{0},y expr=(\thisrowno{2})] {compredunwpghifglim.dat};
     \addplot+[black, mark=triangle*, mark size=1pt, mark options={black}, line width=1pt] table[x expr=\thisrowno{0},y expr=(\thisrowno{1})] {compredunwpghifglim.dat};
\end{axis}
\end{tikzpicture}
\end{subfigure}\caption{
Comparison for different redundancies in WPGHI (top) and W-FGLIM (bottom) for
70 signals sorted by medhigh redundancy performance; the sort sequences are available on the manuscript webpage.\label{fig:redu2}}
\end{figure}

\begin{figure*}[tbp]
\centering
\includegraphics[width=0.6\textwidth,width=0.5\textwidth]{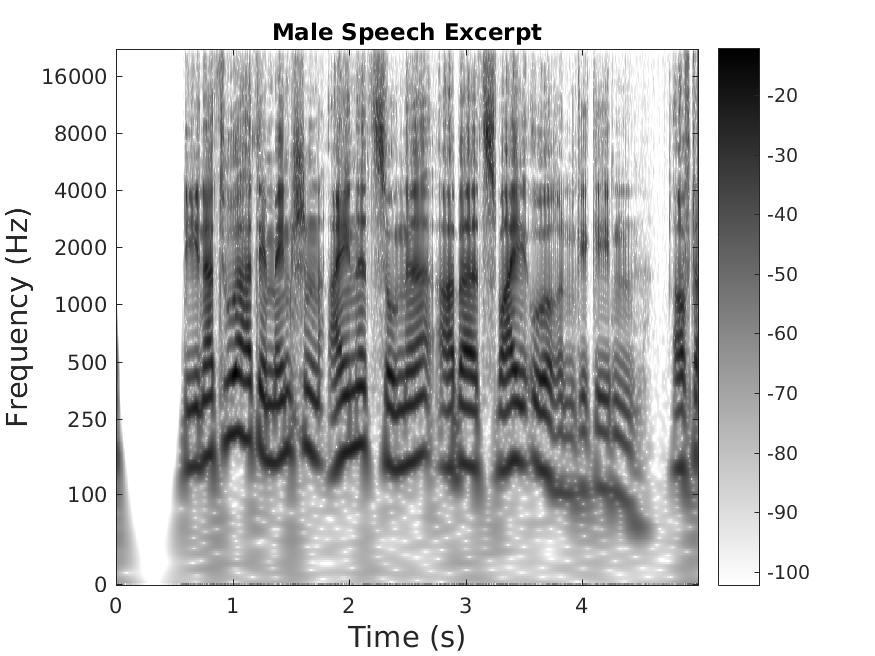}\\
\includegraphics[width=0.6\textwidth,width=0.5\textwidth]{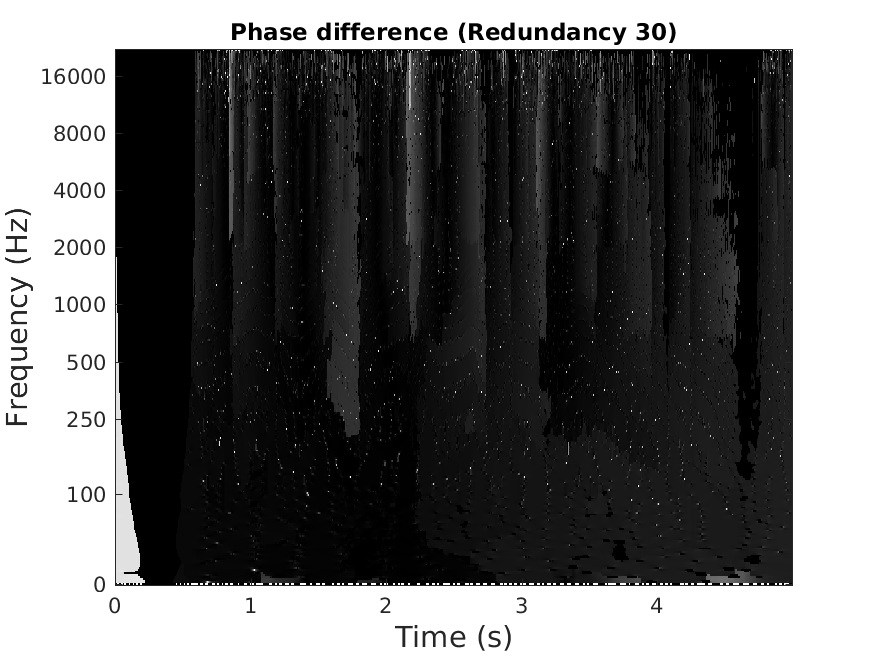}\includegraphics[width=0.6\textwidth,width=0.5\textwidth]{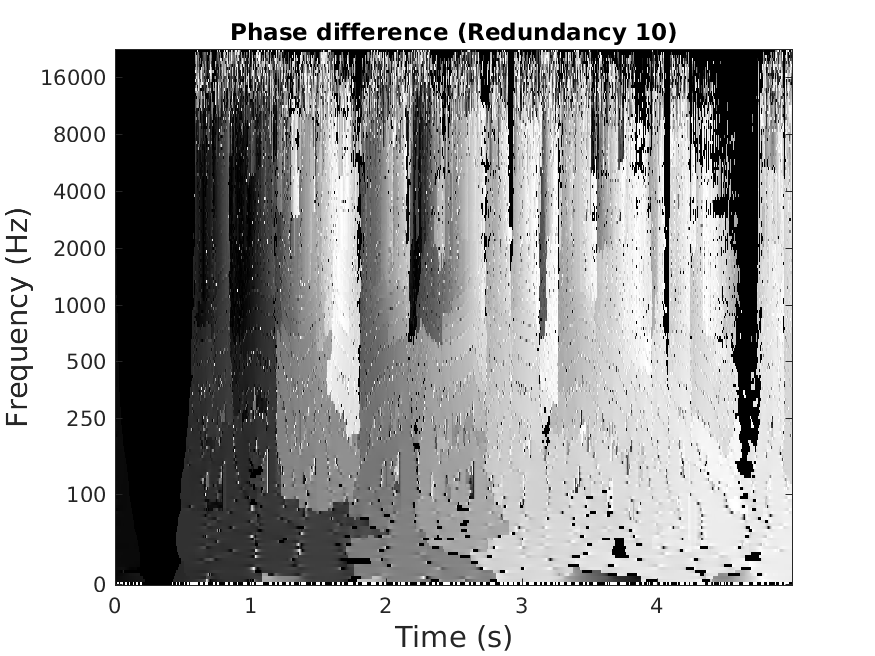}\\
\includegraphics[width=0.6\textwidth,width=0.5\textwidth]{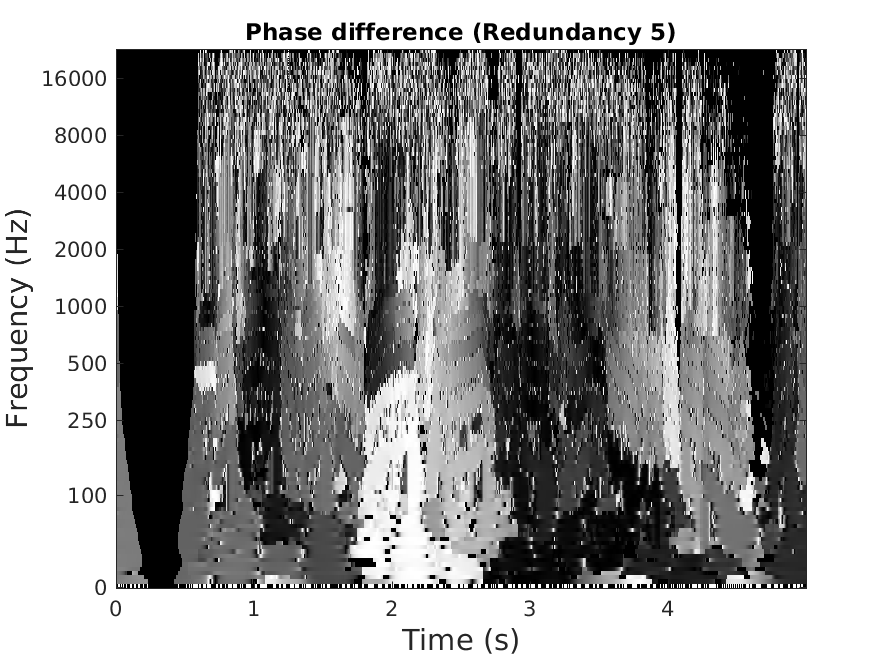}\includegraphics[width=0.6\textwidth,width=0.5\textwidth]{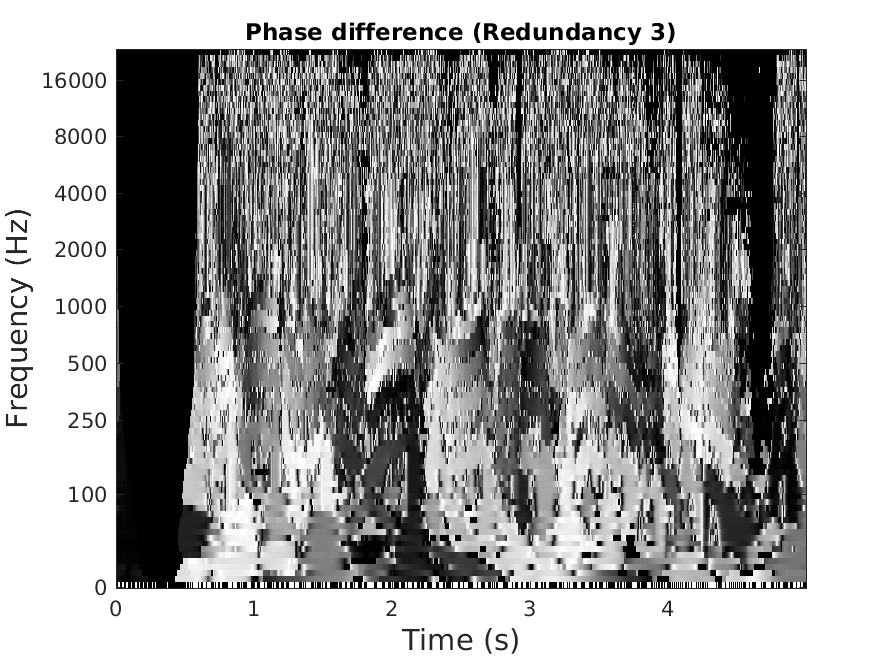}\\
\caption{Top: Wavelet scalogram of a male German speech recording (derived from signal $54$ of the SQAM dataset). Other panels: Difference between the true WT phase the test signal and the phase estimate proposed by WPGHI. Gray level indicates the difference in the range $0$ (black) to $\pi$ (white). Whenever the WT magnitude is below a tolerance level, the phase difference is also set to zero. Phase difference fluctuation clearly increases with decreasing redundancy.\label{fig:phasediffs}}
\end{figure*}

  Similar to Experiment I, mean value and standard deviation over all signals are presented in Table \ref{tab:compRedundancies}, for all parameter sets, with detailed results for all test signals shown in Figure~\ref{fig:redu2}. 
  Additionally, Figure \ref{fig:phasediffs} shows an example of the difference between the target phase and the WPGHI-proposed phase estimate at different redundancies. 
  
  As expected, performance of both proposed methods increases with increasing redundancy.
  These improvements are apparent in the average performance over the whole signal set, but in most cases also on the level of individual signals, as can be seen in Figure~\ref{fig:redu2}. 
  Similar to Experiment I, the average improvement in terms of spectral convergence of W-FGLIM over WPGHI is significant, ranging between $10$ (low) and $6$ dB (high).    
 
 \begin{table}[thb] 
     \begin{center}\begin{tabular}{ |cc|cc|} 
      \hline 
     $K/a_{\mathrm d}$ & method &  WPGHI & W-FGLIM \\ 
      \hline 
     \multirow{2}{*}{30} & mean & $ -38.2012 $ & $ -44.1790 $ \\
      & std & $   6.6897 $ & $   9.8989 $ \\ 
      \hline 
     \multirow{2}{*}{10}& mean & $ -32.5835 $ & $ -40.7522 $ \\
      & std & $   6.5999 $ & $   6.9075 $ \\  
      \hline 
     \multirow{2}{*}{5} & mean & $ -26.8400 $ & $ -36.4934 $ \\ 
      & std & $   5.9723 $ & $   5.6917 $ \\ 
      \hline 
     \multirow{2}{*}{3} &mean & $ -20.7873 $ & $ -29.7372 $ \\ 
      & std & $   5.1843 $ & $   6.6363 $ \\ 
      \hline 
      \end{tabular}\end{center} 
     \caption{Means and standard deviation of spectral convergence for WPGHI and WPGHI-FGLIM for redundancies $K/a = 30,10,5,3$, with $\alpha = 1000$.}
     \label{tab:compRedundancies}
     \end{table} 

     
The performance difference between the different redundancies is also apparent in Figure \ref{fig:phasediffs}. The phase reconstruction quality can be visually estimated from the characteristics of the difference between the target phase and the proposed estimate. Large areas of flat color indicate good quality; as the quality decreases, the phase difference becomes more \emph{patchy}, with stronger fluctuation within patches. The figure shows that this \emph{patchiness} is closely linked to redundancy of the underlying wavelet representation, or more generally, the employed sampling scheme. This behavior is characteristic for WPGHI and was previously observed for STFT-based PGHI~\cite{ltfatnote040} as well.

For all redundancies, informal listening has shown good to excellent reconstruction quality. 
Usually, the reconstructions using medium, medium high, or high redundancy and either of the proposed methods were indistinguishable from one another.
In the case where distortions are audible at all, they were never found to be irritating.  
At low redundancy, the results are generally good as well, but a larger number of examples has clearly audible artifacts, e.g., signals $14$ (oboe), $15$ (cor anglais) and $16$ (clarinet).

 \section{Conclusion}
 
We have highlighted a number of interesting properties of the partial
derivatives of the continuous WT. In particular, we showed
how the log-magnitude and phase partial derivatives are related to wavelet
transforms with modified mother wavelets. 
We characterize the class of wavelets that generate analytic
WTs and obtain as a result a generalization of the Cauchy wavelets. 
Based on the analyticity of these WTs, we obtain explicit phase-magnitude relations.
Similar to the Gaussian STFT, the phase-magnitude
relationship can be used as a basis for implementing magnitude-only
reassignment and phaseless reconstruction. We explored the second
application, providing very good results when applied to complex
audio data, in most cases free of any perceptible distortion. We
demonstrated that reconstruction quality is often notably improved over the
established Griffin-Lim algorithm and our own previous implementation,
relying on approximate phase-magnitude relations for Gaussian filter banks.

Future work will be concerned with real-time implementation of the proposed
algorithm in the style of RTPGHI~\cite{ltfatnote048} and its generalization
to nonuniform decimation schemes. The new phase-magnitude relations for
WTs are crucial for the derivation of more appropriate phase
approximation schemes for more general time-frequency filters, improving and
extending the previously proposed filter bank PGHI algorithm~\cite{ltfatnote051}. 
We also plan to investigate phase-magnitude relations for
the polyanalytic generalizations of the Cauchy WT \cite{Abreu}.
Finally, the proposed scheme can serve as starting point for a
wavelet-based phase vocoder for time-stretching and pitch-shifting of audio
in the spirit of~\cite{ltfatnote050}.
 
  \section*{Acknowledgments}
  We wish to thank Andr\'es Marafioti for performing a small listening test confirming the reported observations on the perceptual quality of the reconstructed audio samples.
  Furthermore, we would like to thank Patrick Flandrin for pointing us to the literature on  Klauder wavelets  which coincide with the analyticity inducing wavelets given by \eqref{eq:solutiondiffeq}.
  Finally, we thank the reviewers, whose comments helped us improve our results and their presentation.

\appendices

\section{Proof of Theorem \ref{thm:derexpr}}\label{app:derivs}
 \begin{proof}
  Under the assumption that $\mwlet\in\LtR$ is continuously differentiable with $\mwlet',\bd T\mwlet'\in\LtR$,
  we can exchange differentiation and integration in \eqref{eq:defwltransform}, see Appendix \ref{app:WTdiff}.
  Thus, the partial derivatives of $W_\mwlet s$ can be expressed as 
  \begin{equation}\label{eq:timederiv1}
  \begin{split}
    \Parx W_\mwlet s(x,y) & = -\frac{1}{y \sqrt{y}} \int_\RR s(t)\overline{\mwlet'\left(\frac{t-x}{y}\right)}\, dt \\
    & = -\frac{1}{y } W_{\mwlet'} s(x,y)
    \end{split}
  \end{equation}
	and 
  \begin{align}\label{eq:scalederiv1}
    & \Pary W_\mwlet s(x,y) \notag \\ 
		& \quad = -\frac{1}{y \sqrt{y}} \int_\RR s(t)\overline{\left(\frac{\mwlet(\bullet)}{2}+(\bullet) \mwlet'(\bullet)\right)\left(\frac{t-x}{y}\right)}\, dt \notag  \\
    & \quad = -\frac{1}{y } \left(\frac{W_\mwlet s(x,y)}{2} + W_{\bd T(\mwlet')} s(x,y)\right) \notag \\
    & \quad = -\frac{1}{y } \left(-\frac{W_\mwlet s(x,y)}{2} + W_{(\bd T \mwlet)'} s(x,y)\right).
  \end{align}
  Here, we used that $(\bd T \mwlet)' = \mwlet+\bd T\mwlet'$ and that that the WT is conjugate linear with respect to the chosen wavelet.
  Using that $\Parx \log(W_\mwlet s) = \frac{\Parx W_\mwlet s}{W_\mwlet s}$ and similarly for the partial derivative with respect to $y$, we obtain for 
  all $(x,y)$ with $W_\mwlet s(x,y) \neq 0$, 
	\begin{subequations}\label{eq:logderiv1}
  \begin{align}
    \Parx \log(W_\mwlet s)(x,y) & = - \frac{1}{y }\frac{W_{\mwlet'} s(x,y)}{W_\mwlet s(x,y)}\label{eq:logtimederiv1}
    \\
    \Pary \log(W_\mwlet s)(x,y) & = \frac{1}{2y} - \frac{1}{y } \frac{W_{(\bd T \mwlet)'} s(x,y)}{W_\mwlet s(x,y)}.\label{eq:logscalederiv1}
  \end{align}
	\end{subequations}
	Taking real and imaginary parts in \eqref{eq:logderiv1} results in \eqref{eq:derlogmag} and \eqref{eq:derphase}, respectively.
\end{proof}

\section{Proof of Theorem~\ref{thm:analyticwlt}}\label{app:analyticwlt}

  We first argue that analyticity of $h$ and differentiability of $f$ already imply that 
  $\mwlet$ satisfies the assumptions in Theorem~\ref{thm:derexpr}, i.e., 
  $\mwlet$ is continuously differentiable with 
	$\mwlet, \mwlet', \bd T\mwlet'\in\LtR$.
  To this end, we note that the assumptions imply that $W_\mwlet s(x-aby,by)$ 
  must be $\mathcal{C}^{\infty}$ in $x$ and $y$ for an arbitrary $s\in \LtR$.
  The same holds true for $W_\mwlet s(x,y)$.
  
  At $x=0$ and $y=1$, the derivative $\Parx W_\mwlet s$
  can be written as 
  \begin{equation}
    \Parx W_{\mwlet}s(0,1)
    = \lim_{x_0\rightarrow 0} \left\langle s, \frac{ \mwlet- \bd T_{ -x_0} \mwlet}{x_0} \right\rangle,
  \end{equation}
  which converges for every fixed $s\in \LtR$. 
  Thus, by a variant of the Banach-Steinhaus theorem \cite[Ch.~II.1, Corollary~2]{yo80}, 
  the limit $\lim_{x_0\rightarrow 0}  \frac{ \mwlet- \bd T_{ -x_0} \mwlet}{x_0}$
  represents a continuous linear functional and hence, by Riesz representation theorem, an element $\mwlet_x$ in $\LtR$.
  Rewriting the derivative for compactly supported, smooth $s$ alternatively as
  \begin{align}
    \Parx W_{\mwlet}s(0,1)
    & = \lim_{x_0\rightarrow 0} 
    \left\langle \frac{ s- \bd T_{ x_0}s}{x_0},  \mwlet \right\rangle
    = -\left\langle  s',   \mwlet \right\rangle,
  \end{align}
  we see that $\mwlet_x$ is the weak derivative of $\mwlet$, i.e., the weak derivative of $\mwlet$ exists and belongs to $\LtR$. Repeating the argument for higher derivatives guarantees that weak derivatives of arbitrary order exist and by standard Sobolev embeddings so do continuous derivatives.
  
  Similarly, the derivative $\Pary W_\mwlet s$ at $x=0$ and $y=1$  
  can be written as
  \begin{equation}
    \Pary W_{\mwlet}s(0,1)
    = \lim_{y_0\rightarrow 0} \left\langle s, \frac{ \mwlet- \bd D_{ -y_0} \mwlet}{y_0} \right\rangle.
  \end{equation}
  Again we have the convergence $\lim_{y_0\rightarrow 0}  \frac{ \mwlet- \bd D_{ -y_0} \mwlet}{y_0} = \mwlet_y$.
  Rewriting the derivative for $s\in \mathcal{C}^{\infty}_{00}$ alternatively as
  \begin{align}
    \Pary W_{\mwlet}s(0,1)
    & = \lim_{y_0\rightarrow 0} 
    \left\langle \frac{ s- \bd D_{{1}/{ -y_0 }}s}{y_0},  \mwlet \right\rangle
    \notag \\
    & = \left\langle s + \bd T s',   \mwlet \right\rangle
    \notag \\
    & = \left\langle (\bd T s)',  \mwlet \right\rangle.
  \end{align}
  Now the operator $s\mapsto  \bd T s'$ is well defined for compactly supported, smooth $s$ with adjoint $s\mapsto  -(\bd T s)'$.
  Furthermore, $\mwlet$ is in the domain of this operator because  $\langle  - (\bd T s)',   \mwlet \rangle =\langle s,  \mwlet_y\rangle$ for all $s$ in a dense subset. 
  Thus, $\bd T \mwlet' = \mwlet_y \in \LtR$.
  Hence, we established all assumptions of Theorem~\ref{thm:derexpr}.

  Analyticity of the function $h$ is equivalent to it satisfying the CR equations that can be compactly expressed as 
    $
  	   \Parx h =-i \Pary h$.
    To rewrite the CR equations for the function $h$ in \eqref{eq:modwlettrans}, we use \eqref{eq:timederiv1}, to obtain
    \begin{align}
      \Parx h 
      & = \bigg(\Parx f (x,y) \bigg) W_\mwlet s(x-aby,by)
      \notag \\*
      &  \quad
      - \frac{ f (x,y)}{by} W_{\mwlet'} s(x-aby,by).
    \end{align}
    Similarly, by \eqref{eq:scalederiv1}, we have
    \begin{align}
      \Pary h
      & = 
      \bigg(\Pary f (x,y)\bigg) W_\mwlet s(x-aby,by)
      \notag \\*
      &  \quad
      - \frac{f (x,y)}{y}  \bigg(-\frac{W_\mwlet s(x-aby,by)}{2} 
      \notag \\* 
      & \quad 
      + W_{(\bd T \mwlet)'} s(x-aby,by) 
      - a W_{\mwlet'} s(x-aby,by)\bigg)
      \notag \\
      & = 
      \bigg( \Pary f (x,y) + \frac{f (x,y)}{2y}\bigg) W_\mwlet s(x-aby,by)
      \notag \\*
      &  \quad
      - \frac{f (x,y)}{y}    W_{(\bd T \mwlet)'} s(x-aby,by) 
      \notag \\* 
      & \quad 
      + \frac{f (x,y)}{y}a W_{\mwlet'} s(x-aby,by)
      .
    \end{align}
    Inserting these expressions into the CR equations results in
    \begin{align}
      & \bigg(\frac{y \Parx f (x,y) +i y \Pary f (x,y)}{f (x,y)} +\frac{i}{2}\bigg)
      W_\mwlet s(x-aby,by)
      \notag \\
      & = 
      \frac{1-iab}{b} W_{\mwlet'} s(x-aby,by)
      + i W_{(\bd T \mwlet)'} s(x-aby,by).
      \label{eq:creqforh}
    \end{align}
    We note that this condition depends on $f$ only via the function
    $g(x,y)= \frac{y \Parx f (x,y) +i y \Pary f (x,y)}{f (x,y)} +\frac{i}{2}$.
    Moreover, using the definition of the WT in \eqref{eq:creqforh} implies that 
    \begin{align}
      \int_{\RR}s(t) 
      \bigg(
      g(x,y) \overline{\mwlet\bigg(\frac{t-x}{y}\bigg)}
      & - \frac{1-iab}{b} \overline{\mwlet'\bigg(\frac{t-x}{y}\bigg)}
      \notag \\
      & \quad
      - i \overline{(\bd T \mwlet)'\bigg(\frac{t-x}{y}\bigg)}
      \bigg)
       dt
       = 
      0
    \end{align}
    for all $s\in \LtR$ and thus
    \begin{align}
      g(x,y) \overline{\mwlet\bigg(\frac{t-x}{y}\bigg)}
      & - \frac{1-iab}{b} \overline{\mwlet'\bigg(\frac{t-x}{y}\bigg)}
      \notag \\
      & \quad
      - i \overline{(\bd T \mwlet)'\bigg(\frac{t-x}{y}\bigg)}
      = 
      0
    \end{align}
    as a function of $t$ in $\LtR$.
    In particular, this implies that $g(x,y)$ must be a constant $w\in \CC$
    and further
    \begin{equation}
      \overline{w} \mwlet -\frac{1+iab}{b} \mwlet' +  i(\bd T \mwlet)'= 0 \,.
      \label{eq:diffeqcauchy}
    \end{equation} 
    To solve this differential equation, it is more convenient to consider the Fourier transformed equivalent of \eqref{eq:diffeqcauchy}. 
		Using the standard properties of the Fourier transform $\widehat{\mwlet'} = 2\pi i \bd T\widehat{\mwlet}$ and $\widehat{\bd T \mwlet} = -(2\pi i)^{-1}(\widehat{\mwlet})'$, this is easily seen to be given by
		\begin{equation}
      \overline{w} \widehat{\mwlet} - 2\pi i \frac{1+iab}{b} \bd T\widehat{\mwlet} -i  \bd T (\widehat{\mwlet})'= 0.
      \label{eq:foudiffeq}
		\end{equation}
		We first note that our assumption $\widehat{\mwlet}(\xi)=0$ for $\xi< 0$ satisfies this differential equation on this domain. 
		For $\xi> 0$ and $\widehat{\mwlet} \neq 0 $, we can easily solve the differential equation by rewriting
		\begin{equation}
		\frac{(\widehat{\mwlet})'(\xi)}{\widehat{\mwlet}(\xi)} 
		= \frac{\overline{w}}{i \xi} - 2\pi \frac{1+iab}{b}
		\label{eq:solvediffeq},
		\end{equation}
which gives
		\begin{align}
		\widehat{\mwlet}(\xi) 
		& = c e^{- i \overline{w} \log  \xi} e^{ - 2\pi \frac{1+iab}{b} \xi}
    \notag \\
		& = c \xi^{-\operatorname{Im}(w)} e^{- 2\pi (\frac{1}{b}+ia) \xi} e^{- i \operatorname{Re}(w) \log  \xi}
		\end{align}
    %
    %
    %
		for an arbitrary constant $c\in \CC$.
    To obtain the parameters used in the theorem, we substitute 
    $- \operatorname{Im}(w) = \frac{\alpha-1}{2}$, $\frac{1}{b}= \regamma$,  $a = \imgamma$, and $\operatorname{Re}(w) = -\beta$.
    Based on our assumptions, we obtain the constraints $\regamma >0$ and $\alpha > -1$ to guarantee $\widehat{\mwlet}\in \LtR$. 

\section{Proof of Theorem~\ref{thm:pmrel}}\label{app:PMrels}

	We can use the CR equations to obtain relationships between the derivatives of real and imaginary parts of the WT, which we will show to result in \eqref{eq:xDirPhaseDeriv} and \eqref{eq:yDirPhaseDeriv}.
	For an arbitrary analytic function $h=u+iv$ the CR equations hold and are given by
	$
	\Parx u = \Pary v 
	$
	and
	$
	\Pary u =- \Parx v$.
	Writing $h=M e^{i\phi}$, the CR equations imply that
	\begin{align}
	 \Parx \phi
	& = - \Pary \log M
	\label{eq:CRphix}
	\end{align}
	and 
	\begin{equation}
	\Pary \phi
	 = \Parx \log M.
	\label{eq:CRphiy}
	\end{equation}
	Using \eqref{eq:CRphix} and \eqref{eq:CRphiy} for the function given by \eqref{eq:genanalytic}, yields
  \begin{align}
    & \Parx \bigg(\phi_{\mwlet}^s\bigg(x- \frac{\imgamma }{\regamma } y,\frac{y}{\regamma }\bigg) + \beta \log y\bigg)
    \notag \\*
    & = - \Pary  \log \bigg(y^{-\frac{\alpha}{2}}\magmwlet\bigg(x- \frac{\imgamma }{\regamma } y,\frac{y}{\regamma } \bigg)\bigg) 
  \end{align}
  and
 	\begin{align}
    & \Pary \bigg( \phi_{\mwlet}^s \bigg(x- \frac{\imgamma }{\regamma } y,\frac{y}{\regamma }\bigg) + \beta \log y \bigg)
    \notag \\
    & = \Parx \log \bigg(y^{-\frac{\alpha}{2}}\magmwlet\bigg(x- \frac{\imgamma }{\regamma } y,\frac{y}{\regamma } \bigg) \bigg).
	\end{align}
  These are equivalent to
  \begin{align}
	  \Parx  \phi_{\mwlet}^s
    & =   \frac{\alpha}{2  y \regamma }
    - \frac{1}{\regamma } \Pary  \log \big(\magmwlet\big) 
    \notag \\*
    & \quad 
    + \frac{\imgamma }{\regamma } \Parx  \log \big(\magmwlet\big)
  \label{eq:phiparxexpr}
	\end{align}
	and
  \begin{equation}
    - \frac{\imgamma }{\regamma }\Parx  \phi_{\mwlet}^s
    + \frac{1}{\regamma } \Pary  \phi_{\mwlet}^s  + \frac{\beta}{ y\regamma }
    = \Parx  \log  \big(\magmwlet  \big).
    \label{eq:phiparyexpr}
	\end{equation}
  Inserting \eqref{eq:phiparxexpr} into \eqref{eq:phiparyexpr} finally results in
  \begin{align}
    & - \frac{\imgamma }{\regamma } \bigg(\frac{\alpha}{2  y \regamma }
    - \frac{1}{\regamma } \Pary  \log \big(\magmwlet\big) 
    \notag \\*
    & 
    + \frac{\imgamma }{\regamma } \Parx  \log \big(\magmwlet\big)\bigg)
    + \frac{1}{\regamma } \Pary  \phi_{\mwlet}^s  + \frac{\beta}{ y\regamma }
    \notag \\
    & = \Parx  \log  \big(\magmwlet  \big),
  \end{align}
  which is equivalent to
  \begin{align}
    \regamma  \Pary  \phi_{\mwlet}^s 
    & = 
    \frac{\alpha\imgamma }{2  y }  - \frac{\beta}{ y }
    + \lvert \mygamma\rvert^2  \Parx  \log  \big(\magmwlet \big)
    \notag \\*
    & \quad 
    - \imgamma     \Pary  \log \big(\magmwlet\big) .
  \end{align}
  This concludes the proof.

  \section{On Differentiability of the Wavelet Transform}\label{app:WTdiff}
  
  Denote by $\Phi_z$ and $\Psi_z$ the difference operators
  \[
   \Phi_z \psi = \frac{\bd T_z \psi - \psi}{z}\ \text{ and }\ \Psi_z \psi = \frac{\bd D_{1+(z-1)} \psi - \psi}{z-1},
  \]
  and note that 
  \[
   \begin{split}
     \frac{\partial}{\partial x} W_{\mwlet} s(x,y) 
   & = \lim_{x_0\rightarrow 0} \frac{W_{\mwlet} s(x+x_0,y)-W_{\mwlet} s(x,y)}{x_0}\\
   & = \lim_{x_0\rightarrow 0} \frac{1}{y}\left\langle s, \bd T_x \bd D_y (\Phi_{x_0/y}\mwlet)\right\rangle,
   \end{split}
  \]
  as well as 
  \[
   \begin{split}
     \frac{\partial}{\partial y} W_{\mwlet} s(x,y) 
   & = \lim_{y_0\rightarrow 0} \frac{W_{\mwlet} s(x,y+y_0)-W_{\mwlet} s(x,y)}{y_0}\\
   & = \lim_{y_0\rightarrow 0} \frac{1}{y}\left\langle s, \bd T_x \bd D_y (\Psi_{1+y_0/y}\mwlet)\right\rangle,
   \end{split}
  \]
  provided that the right-hand sides converge.
  
  We first show that, for $z \rightarrow 0$
  \begin{equation}
   \Phi_{z}\mwlet \rightarrow -\mwlet' \text{ and } \Psi_{1+z}\mwlet \rightarrow -(\mwlet/2 + \bd T \mwlet')
  \label{eq:convergencephipsi}
  \end{equation}
  as functions in $\bd L^2(\RR)$, provided that $\mwlet$, $\mwlet'$, and $\bd T \mwlet'$ are elements of 
  $\bd L^2(\RR)$.

  Since $\mwlet'$ is assumed to be in $\bd L^2(\RR)$, there is for every $\tilde{\epsilon}>0$ an $r_{\tilde{\epsilon}}>0$,
  such that $\big\lVert\mwlet'\chi_{\RR\setminus \overline{B_{r_{\tilde{\epsilon}}}(0)}}\big\rVert < \tilde{\epsilon}$.
  Moreover, by the fundamental theorem of calculus,
  \begin{align} 
    |\Phi_{z}\mwlet(t) + \mwlet'(t)| 
    & = \bigg\lvert \frac{\int_t^{t-z} \mwlet'(s)\,ds}{z}   + \mwlet'(t) \bigg\rvert
    \notag \\
    & = \bigg\lvert - \int_0^1 \mwlet'(t-sz) + \mwlet'(t) \,ds \bigg\rvert
    \label{eq:phizppsip} 
  \end{align}
  Using \eqref{eq:phizppsip}, we obtain
  \begin{align}
   & \Big\lVert (\Phi_{z}\mwlet + \mwlet') \chi_{\RR\setminus \overline{B_{r_{\tilde{\epsilon}}+|z|}(0)}}\Big\rVert^2 
  \notag \\ 
  & \leq 
  \int_{\RR\setminus \overline{B_{r_{\tilde{\epsilon}}+|z|}(0)}} \bigg( \int_0^1 \lvert \mwlet'(t-sz) - \mwlet'(t) \rvert \,ds \bigg)^2 dt
  \notag \\ 
  & \leq 
  \int_{\RR\setminus \overline{B_{r_{\tilde{\epsilon}}+|z|}(0)}}  \int_0^1 3\lvert \mwlet'(t-sz)\rvert^2 + 3 \lvert \mwlet'(t) \rvert^2 \,ds  \, dt
  \notag \\ 
  & < 
  6 \tilde{\epsilon}^2
  \label{eq:boundouter}
  \end{align}
  where we used Jensen's inequality and Fubini's theorem.
  Furthermore, we have $|\mwlet'(t-s)-\mwlet'(t)|<\epsilon_{r}(|s|)$ for all $|t|<r$ and $|s|<1$ by uniform continuity of $\mwlet'$ on the compact set $B_{r +1}(0)$, where $\epsilon_r( \delta)\searrow 0$ for $\delta\rightarrow 0$.
  Thus, similar to \eqref{eq:boundouter}, we obtain
  \begin{align}
    \Big\lVert  (\Phi_{z}\mwlet + \mwlet') \chi_{\overline{B_{r_{\tilde{\epsilon}}+|z|}(0)}}\Big\rVert^2  
    & <  2(r_{\tilde{\epsilon}}+|z|) \epsilon_{r_{\tilde{\epsilon}}+1}^2( |z|)
  \end{align}
  for $|z|<1$.
  Finally, for every $\tilde{\epsilon}$, there is a $z_{\tilde{\epsilon}}\in (0,1)$, such that $\epsilon_{r_{\tilde{\epsilon}}+1}(z_{\tilde{\epsilon}}) < \tilde{\epsilon}/\sqrt{2(r_{\tilde{\epsilon}}+1)}$, 
  which implies $\| \Phi_{ z}\mwlet + \mwlet'\| <  \sqrt{7}\tilde{\epsilon}$ for all $|z|\leq z_{\tilde{\epsilon}}$.
  
  A similar, but slightly more complicated argument shows that $\Psi_{1+z}\mwlet \rightarrow -(\mwlet/2 + \bd T \mwlet')$, 
  provided $\mwlet, \bd T \mwlet'\in\bd L^2(\RR)$.
  Here, we start with an $r_{\tilde{\epsilon}}>0$,
  such that $\| \bd T \mwlet'\chi_{\RR\setminus \overline{B_{r_{\tilde{\epsilon}}}(0)}}\| < \tilde{\epsilon}$ and $\|  \mwlet \chi_{\RR\setminus \overline{B_{r_{\tilde{\epsilon}}}(0)}}\| < \tilde{\epsilon}$.
  As above, we use the fundamental theorem of calculus to obtain
  \begin{align} 
    & |\Psi_{1+z}\mwlet(t) + \mwlet(t)/2 + \bd T \mwlet'(t)| 
    \notag \\
    & = \Bigg\lvert \frac{\frac{1}{(1+z)^{1/2}}\mwlet\big(\frac{t}{1+z}\big)-\mwlet(t)}{z}
      +\frac{\mwlet(t)}{2} + t \mwlet'(t) \Bigg\rvert
    \notag \\
    & = \Bigg\lvert \frac{\int_0^{z} \Big[ \frac{1}{(1+\bullet)^{1/2}}\mwlet\big(\frac{t}{1+\bullet}\big)\Big]'(s) \,ds}{z}   
      +\frac{\mwlet(t)}{2} + t \mwlet'(t) \Bigg\rvert
    \notag \\
    & = \Bigg\lvert - \int_0^{z} \frac{ \frac{1}{2}\mwlet\big( \frac{t}{1+s} \big) + \frac{t}{1+s} \mwlet'\big( \frac{t}{1+s} \big)  }{(1+s)^{3/2}z} \,ds  
      +\frac{\mwlet(t)}{2} + t \mwlet'(t) \Bigg\rvert
    \notag \\
    & = \Bigg\lvert - \int_0^{1} \frac{ \frac{1}{2}\mwlet\big( \frac{t}{1+sz} \big) + \frac{t}{1+sz} \mwlet'\big( \frac{t}{1+sz} \big)  }{(1+sz)^{3/2}}
      +\frac{\mwlet(t)}{2} + t \mwlet'(t)  \,ds  \Bigg\rvert
    \notag \\
    & \leq   \int_0^{1} \Bigg\lvert \frac{ \frac{1}{2}\mwlet\big( \frac{t}{1+sz} \big)  }{(1+sz)^{3/2}} - \frac{\mwlet(t)}{2}  \Bigg\rvert +  \Bigg\lvert \frac{ \frac{t}{1+sz} \mwlet'\big( \frac{t}{1+sz} \big)  }{(1+sz)^{3/2}}
       - t \mwlet'(t)  \Bigg\rvert \,ds 
    \label{eq:boundpsiabs}
  \end{align}
  Similar to \eqref{eq:boundouter}, this results in 
  \[
   \Big \lVert (\Psi_{1+z}\mwlet + \mwlet/2 + \bd T\mwlet' )\chi_{\RR\setminus \overline{B_{(1+|z|)r_{\tilde{\epsilon}}}(0)}} \Big \rVert^2 < 16\tilde{\epsilon}^2 .
  \]
  Furthermore, we can bound for $t\in \overline{B_{2r}(0)}$
  \begin{align}
    & \Bigg\lvert \frac{ \frac{1}{2}\mwlet\big( \frac{t}{1+s} \big)  }{(1+s)^{3/2}} - \frac{\mwlet(t)}{2}  \Bigg\rvert
    \notag \\
    & \leq 
     \frac{1}{2 \lvert(1+s)^{3/2}\rvert}\Bigg\lvert \mwlet\bigg( \frac{t}{1+s} \bigg)  -  \mwlet ( t)    \Bigg\rvert + \frac{\lvert \mwlet ( t) \rvert }{2} \Bigg\lvert  \frac{1}{(1+s)^{3/2}} - 1 \Bigg\rvert
    \notag \\
    & \leq  
    \epsilon_r(t  \lvert s\rvert ) 
    + \frac{\lvert \mwlet ( t) \rvert }{2} \epsilon_r( \lvert s\rvert )
  \end{align}
  for any $\lvert s\rvert\leq 1/2$, where $\epsilon_r (\delta) \to 0$ monotonically  for $\delta \to 0$.
  Analogously, we obtain
  \begin{align}
     \Bigg\lvert \frac{ \frac{t}{1+s} \mwlet'\big( \frac{t}{1+s} \big)  }{(1+s)^{3/2}}
       - t \mwlet'(t)  \Bigg\rvert
    & \leq 
    2 \epsilon_r(t  \lvert s\rvert ) + \lvert \mwlet ( t) \rvert   \epsilon_r( \lvert s\rvert ).
    \label{eq:part2boundpsiabs}
  \end{align}
  Equations \eqref{eq:boundpsiabs}--\eqref{eq:part2boundpsiabs} imply 
  \begin{equation} 
     \bigg\lvert\Psi_{1+z}\mwlet(t) + \frac{\mwlet(t)}{2} + \bd T \mwlet'(t)\bigg\rvert 
    \leq  3 \epsilon_r(t  \lvert z\rvert ) + 2 \lvert \mwlet ( t) \rvert  \epsilon_r( \lvert z\rvert ) 
  \end{equation}
  and thus
  \begin{align}
    & \Big\lVert \Big(\Psi_{1+z}\mwlet +  \frac{\mwlet}{2} + \bd T\mwlet'\Big) \chi_{\overline{B_{(1+|z|)r_{\tilde{\epsilon}}}(0)}} \Big\rVert 
    \notag \\ 
    & < \sqrt{4 r_{\tilde{\epsilon}}}
    \bigg(3+ 2 \sup_{t\in \overline{B_{2 r_{\tilde{\epsilon}}}(0)}} \lvert \mwlet(t) \rvert \bigg)\epsilon_{r_{\tilde{\epsilon}}}\big(2 r_{\tilde{\epsilon}}\lvert z\rvert\big).
  \end{align}
  for $\lvert z\rvert<1/(4r_{\tilde{\epsilon}})$ and where we assumed for simplicity $2r_{\tilde{\epsilon}}>1$.
  Again, choosing $\lvert z \rvert$ sufficiently small implies the proposed convergence.
  
  Thus, we finished the proof of \eqref{eq:convergencephipsi}, which implies  
  \[
   \left\langle s, \lim_{x_0\rightarrow 0} \bd T_x \bd D_y (\Phi_{x_0/y}\mwlet)\right\rangle = -\left\langle s, \bd T_x \bd D_y \mwlet'\right\rangle
  \]
  and
  \[
   \left\langle s, \lim_{y_0\rightarrow 0} \bd T_x \bd D_y (\Psi_{1+y_0/y}\mwlet)\right\rangle = -\left\langle s, \bd T_x \bd D_y (\mwlet/2 + \bd T \mwlet')\right\rangle.
  \]
  Because convergence is in  $\bd L^2(\RR)$, we can exchange limit and integral by continuity of the inner product.
  
  Clearly, this argument can be repeated to obtain higher order derivatives, provided that $\mwlet$ has sufficient regularity and decay.

\bibliographystyle{IEEEtran}
\bibliography{addbib}

\end{document}